
%
%
%
%
%
%
%
\tolerance=10000
\magnification=1200
\input amstex
\documentstyle{amsppt}
%
\hoffset=-0.4truein
\voffset=0.2truein
\hsize=7.5truein
\vsize=9.0truein
\addto\tenpoint{\normalbaselineskip=18pt\normalbaselines}
\addto\eightpoint{\normalbaselineskip=15pt\normalbaselines}
%
%
\def\Bl{1}
\def\Bo{2}
\def\Coa{3}
\def\Cob{4}
\def\Con{5}
\def\Da{6}
\def\ENNa{7}
\def\ENNb{8}
\def\Ns{9}
\def\Ph{10}
\def\RS{11}
\def\Scha{12}
\def\Schb{13}
\def\Tr{14}
%
%
\def\varf{f}
\def\K{\Cal K}
\def\pa{\parallel}
\def\Cinfa{C_{\infty}A}
\def\C{\Bbb C}
\def\R{\Bbb R}
\def\N{\Bbb N}
\def\Z{\Bbb Z}
\def\T{\Bbb T}
\def\Cinfc{C_{\infty}\C} %
\def\Cinfap{I_{\infty}A}
\def\Cinfcp{I_{\infty}\C} %
\def\schr{\Cal S(\R)}
\def\sona{\Cal S(\R, A, \alpha)}
\def\Ialph{I_{\infty}(\alpha)} %
\def\Calph{C_{\infty}(\alpha)}

\def\Kinf{\Cal K^{\infty}}
\def\alphh{{\widehat{\alpha}}}

\def\w{\omega}
\refstyle{A}  %
\rightheadtext{Representable $K$-theory of Smooth}
\leftheadtext{N.C. Phillips and L.B. Schweitzer}
\topmatter
\title Representable $K$-theory of Smooth Crossed
Products by $\R$ and $\Z$ \endtitle
\author N. Christopher $\text{Phillips}^{1}$
and Larry B. Schweitzer \endauthor
\address Department of Mathematics, University of
Oregon, Eugene, OR 97403-1222 USA \endaddress
\email phillips\@ bright.uoregon.edu
\endemail
\thanks 1) Research partially supported by NSF
grant DMS 91-06285. \endthanks
\address Department of Mathematics, University of Victoria,
Victoria, B.C. Canada V8W 3P4 \endaddress
\email lschweit\@ sol.uvic.ca
\endemail
\toc
\widestnumber\head{2.20}
\widestnumber\specialhead{3.0}
\specialhead Introduction \endspecialhead
\specialhead 1.  Smooth Crossed Products by $\R$ \endspecialhead
\head 1.1 Definitions and Preliminaries \endhead
%
%
\head 1.2  The Thom Map $\theta$ is an Isomorphism \endhead
\head 1.3  Replacing $A^{\infty}$ with $A$ and $\schr$ with
   $L_{1}(\R)$ \endhead
\specialhead 2.  Smooth Crossed Products by $\Z$ \endspecialhead
\specialhead 3.  Examples and Applications \endspecialhead
\endtoc
\endtopmatter
\document
\heading Abstract \endheading
\par
We show that the Thom isomorphism and the Pimsner-Voiculescu
exact sequence both hold for smooth crossed products of
Fr\'echet algebras by $\R$ and $\Z$ respectively.
We also obtain the same results for $L^{1}$-crossed
products of Banach algebras by $\R$ and $\Z$.
\vskip\baselineskip
\heading Introduction \endheading
\par
In \cite{\Coa}, Connes proved that the
$K$-theory $K_{*}(C^{*}(\R,A))$  is isomorphic to $K_{*+1}(A)$ for any
continuous action of $\R$ on a C*-algebra $A$. In this paper, we prove
the analog of this result for the representable $K$-theory of smooth
crossed products for actions of $\R$ on Fr\'echet algebras. In order
that the representable $K$-theory of the Fr\'echet algebra $A$ be
defined, we assume that all our Fr\'echet algebras are locally
$m$-convex. In order to ensure that the crossed
products are again locally $m$-convex, we assume that our actions are
$m$-tempered (Definition 1.1.3).
No other restrictions are necessary. As a corollary, we obtain
the Pimsner-Voiculescu exact sequence for actions of $\Z$ on Fr\'echet
algebras. As further corollaries, we compute the K-theory of the
$L^1$ crossed products for isometric actions of $\R$
and $\Z$ on arbitrary Banach algebras.
\par
Our smooth crossed products $\Cal S(\R, A)$ and $\Cal S(\Z, A)$
are the standard sets of $A$-valued
Schwartz functions on $\R$ and $\Z$ respectively, with convolution
multiplication.  (In the language of
\cite{\Scha} and  \cite{\Schb}, these are just the smooth
crossed products of Schwartz functions which vanish rapidly with
respect to the gauge $\sigma(x) = |x|$ on $\Z$ or $\R$.
In a future paper, we hope to consider the $K$-theory of
smooth crossed products defined using stronger
rapid vanishing conditions.)
%
%
\par
Fr\'echet algebras occuring as dense subalgebras of
C*-algebras have been important in many recent results in C*-algebra
theory.  For example, the  theory of  differential geometry on a
noncommutative space \cite{\Cob} requires the use of
\lq\lq differentiable structures\rq\rq\
for these noncommutative spaces,
or some sort of algebra of \lq\lq differentiable functions\rq\rq\
on the noncommutative space.  Such algebras of functions have usually
been provided by a dense Fr\'echet subalgebra of smooth functions $A$
for which the $K$-theory $K_{*}(A)$ is the same as the
$K$-theory of the C*-algebra $K_{*}(B)$ (see for example  \cite{\Bl},
\cite{\Bo}, \cite{\ENNa}, \cite{\Ns},
the recent works of  Baum and Connes, Blackadar and Cuntz,
R. Ji, P. Jolissaint, V. Nistor and many others).
In \cite{\Ns} and \cite{\ENNa}, methods are given to compute
the cyclic cohomology \cite{\Cob} of
smooth crossed products of Fr\'echet algebras by $\R$ and $\Z$.
\par
If the Fr\'echet algebras are spectral invariant dense subalgebras of
C*-algebras, and the smooth crossed products are spectral invariant in
the
corresponding C* crossed products, then our results actually follow from
the corresponding C*-algebra results, since the inclusions of the
dense subalgebras into their respective C*-algebras are all
isomorphisms on $K$-theory. However, in \cite{\Ns} no spectral
invariance
assumptions are made on the algebra, and in \cite{\ENNa} the Fr\'echet
algebra need not be given as a dense subalgebra of a Banach algebra.
We also make none of these assumptions.  We
show that the $K$-theory of the smooth crossed product is always related
to that of the original algebra in the expected way, under essentially
the minimum assumptions needed to make sure that all
the Schwartz crossed products are $m$-convex Fr\'echet algebras,
so that a well behaved $K$-theory can be defined for them.
\par
Our proof of the smooth Connes isomorphism is based on the method of
\cite{\ENNb}.  As a proof of the Connes isomorphism which assumes Bott
periodicity, our proof actually reduces in the C*-algebra case to an
even simpler proof than the one in \cite{\ENNb} , because the continuous
fields have been replaced
by ordinary short exact sequences. For more details on the relation with
\cite{\ENNb}, see the beginning of Section 1.
\par
%
%
In \S 1.1, we define a map from the representable $K$-theory
$RK_{*}(A)$ of $A$ to
$RK_{*+1}(\Cal S(\R, A, \alpha))$,
where $\alpha$ is a differentiable action of $\R$ on $A$.  In \S 1.2,
we show that this map is an isomorphism.  In \S 1.3, we show that
the hypothesis that $\alpha$ be a differentiable action
can be dropped, as long as $\alpha$ is assumed strongly
continuous.
We also prove the results for the $L^1$ crossed product, using spectral
invariance arguments.
\par
In \S 2, we define the smooth mapping torus and use this to
show that the Pimsner-Voiculescu exact sequence holds for
smooth crossed products by
$\Z$, following the method of Theorem 10.2.1 of \cite{\Bl}.
In \S 3, we give some applications and examples,
including some nonstandard ones.
\par
We will use the notation $\R$, $\C$,  $\N$, $\Z$, $\T$  for
the reals, complexes, natural numbers (including zero), integers,
and the circle group respectively.
The spaces $\schr$ and $\Cal S (\Z)$ will denote the standard
Fr\'echet space of Schwartz functions on $\R$ and $\Z$ respectively.
All Fr\'echet spaces will be assumed locally convex,
and to simplify terminology we will call a locally $m$-convex
Fr\'echet algebra a Fr\'echet algebra.  ($m$-convex means that
the topology can be given by a family of submultiplicative
seminorms.)  The tensor product $\otimes$
will always mean the completed projective tensor product of
two Fr\'echet spaces.  By an action $\alpha$ of a group $G$ on a
Fr\'echet algebra $A$ we will mean a group homomorphism
$\alpha \colon G \rightarrow {\text{Aut}}(A)$.
The notation $A^{+}$ will mean \lq\lq $A$ with unit adjoined\rq\rq\ ,
even if the algebra $A$ already has a unit.
\par
When we say \lq\lq $K$-theory \rq\rq\  or \lq\lq isomorphism
on $K$-theory \rq\rq, we will mean with regard to the
representable $K$-functor (denoted by $RK$) of \cite{\Ph}.
\par
The first author would like to thank Toshikazu Natsume for supplying him
with a preprint of \cite{\ENNb},
and Theodore Palmer for suggesting some examples.
The second author would like to thank
Theodore Palmer and the first author for a pleasant stay at
the University of Oregon.
%
%
%
%
\heading \S 1 Smooth Crossed Products by $\R$ \endheading
\par
In this section we prove that the Thom isomorphism holds for
smooth crossed products $\sona$, where $\alpha$ is any
strongly continuous (not necessarily smooth), $m$-tempered
action of $\R$ on a Fr\'echet algebra $A$.
\par
Before we begin \S 1.1,
we give a brief discussion of the development of our method of proof.
In \cite{\ENNb}, a map from the suspension of $A$ to the C*-crossed
product is defined roughly as follows.  For $ s \in [0, 1]$,
let $C^{*}(\R,A,  \alpha_s)$ be the C*-crossed product of $A$ with
$\R$, with action $t \mapsto \alpha_{ts}$.
For $s=0$, this is just the suspension of $A$. These algebras fit
together to form a continuous field of C*-algebras
over $[0, 1]$. Let $u$ be some invertible in the unitization of
the suspension of $A$. Find a continuous section $s \mapsto u_s$ through
$u$. Then for $s$ sufficiently small, $u_s$ will  be invertible
in $C^{*}(\R,A, \alpha_s)^{+}$.  The argument for this
depends on the fact that a Banach algebra
has an open invertible group.  Then one notes that for $s>0$, the fibers
$C^{*}(\R,A,\alpha_s)$ are all isomorphic as algebras, so that $u_s$ in
fact defines an invertible in $C^{*}(\R, A,\alpha)^{+}$.
This then defines a map from $K_1(SA)$ to $K_1(C^{*}(\R,A,\alpha))$.
This map is composed with the Bott map to get the Thom map from
$K_0(A)$ to $K_1(C^{*}(\R,A,\alpha))$.
\par
The problem with imitating this for smooth crossed products
is that $\Cal S(\R,A,\alpha_s)^{+}$  may {\it not} have an open group
of invertibles.  So we would not know that for small $s$
the element $u$ above is invertible in $\Cal S(\R,A, \alpha_s)^{+}$.
Moreover, the decomposition of $\Cal S(\R,A,\alpha_s)$
as an inverse limit of Banach algebras may vary as $s$ varies,
making it awkward to work on the Banach algebra level.
We use the following device to overcome this.  We let
$\R$ act on the pointwise multiplication algebra $C^{\infty}([0,1], A)$
by $\beta_t(f)(s) = \alpha_{st}(f(s))$.
Then $\Cal S(\R, C^{\infty}([0, 1], A), \beta)$
is a Fr\'echet algebra which serves as a substitute for the
continuous field $\{\Cal S(\R,A, \alpha_{ts})\}_{s\in [0,1]}$.
There are natural maps from $\Cal S(\R, C^{\infty}([0, 1], A), \beta) $
to the suspension $\Cal S(\R) \otimes A$ and to $\Cal S(\R, A, \alpha)$
given by evaluation at $0$ and $1$ respectively.  We show that the
former map can be inverted in Corollary 1.1.12 below.  Thus we define
our Thom map for smooth crossed products to be this inverse
composed with the second map of evaluation at $1$.
\par
Our proof in \S 1.2 that this composition is an isomorphism proceeds
much like the proof of \cite{\ENNb}, using commutative diagrams and
showing that one edge of the diagram is essentially the identity map,
so that another path of the diagram (which has the Thom map as a
composite) must be also.
\par
We conclude this discussion by looking at the possibility of imitating
other proofs of the Thom isomorphism for C*-crossed products by $\R$,
in particular \cite{\Coa} and \cite{\Bl, \S 10.9}. One essential part of
the proof in \cite{\Coa} is as follows.  One has the crossed product
$C^{*}(\R, A, \alpha)$, and an idempotent $e$ in $A$.  One then wants to
define
an exterior equivalent action $\gamma$ of $\R$ on $A$, such that the
crossed products $C^{*}(\R, A,\alpha)$ and $C^{*}(\R,A,  \gamma)$
are isomorphic, and $\gamma$ leaves $e$ fixed.
This exterior equivalence is implemented by a multiplier
on $C^{*}(\R, A, \alpha)$ given by the formula
$$ u_{t} =
 \sum_{n=0}^{\infty} i^{n} \int_{0\leq s_{1}\leq\dots s_{n} \leq t}
\alpha_{s_{1}}(x) \dots \alpha_{s_{n}}(x) ds_{1}\dots ds_{n},\tag *$$
where $x$ is some self-adjoint element of $A$ depending on $e$.
This formula for $u_{t}$ converges in the norm on $A$, essentially
because
the series $\sum_{n=0}^{\infty} t^{n}/n!$ converges for any $t$.
Also, $u_{t}$ is a unitary element of the C*-algebra $A$
for each $t$. (See proof of Proposition 4 of \cite{\Coa, \S II}.)
However, if we try to imitate this in our setting,
when $A$ is a Banach or Fr\'echet algebra, and we use the
smooth crossed product in place of the C*-crossed product,
then it is not clear that $(*)$ gives us a nice enough $u_{t}$.
In particular we need $\pa u_{t} \pa $ to be bounded by some
polynomial in $t$, whereas the formula $(*)$  seems to
give an exponential bound at best.
\par
Imitating Blackadar's proof \cite{\Bl, \S 10.9} seems to require an
explicit formula for a smooth version of the function $F$ in \cite{\Bl,
Lemma 10.9.3}.  The spectral invariance of the smooth compact operators
$\Kinf$ in the compact operators $\K$ should at least come
very close to ensuring the existence of such an $F$, but
an explicit formula  seems to be difficult to find.
\vskip\baselineskip
\heading \S 1.1 Definitions and Preliminaries  \endheading
\par
The main purpose of this section is to construct a
map $\epsilon_{0}$ of Fr\'echet algebras associated with a smooth action
of $\R$, and to prove that it is an isomorphism on $K$-theory.
This is necessary because $(\epsilon_{0})^{-1}_{*}$
is used in the next section in our definition of the Thom map from
the $K$-theory of $A$ to the $K$-theory of the smooth crossed product.
\subheading{1.1.1 Definition}
Let $A$ be a Fr\'echet algebra. Then the {\it smooth interval
algebra} $\Cinfap$ of $A$ is the Fr\'echet algebra
$C^{\infty}([0, 1], A)$ of all $C^{\infty}$ functions $\varf$ from
$[0, 1]$ to $A$ with the topology of uniform convergence of each
derivative of $\varf$ in each continuous seminorm on $A$.
The {\it smooth cone} of $A$ is
$$\Cinfa = \{\, \varf \in \Cinfap\,
|\,  \varf^{(n)}(0) = 0\, {\text{ for all }} n\in \N \,\},$$
where $\varf^{(n)}$ denotes the $n$th derivative of $\varf$.
We give $\Cinfa$ the topology inherited from $\Cinfap$.
\proclaim{1.1.2 Lemma}
\roster
\item
For each $k, n \in \N$, the seminorm
$$ \pa \varf \pa_{k, n, m}^{\prime} = \sup_{r \in [0, 1]}
\pa r^{-k} \varf^{(n)}(r) \pa_{m},  $$
where $\pa \quad \pa_{m}$ is any continuous seminorm on $A$,
is well-defined and continuous on $\Cinfa$.
\item
The Fr\'echet spaces $\Cinfc$ and $\Cinfcp$ are
nuclear topological vector spaces.
\item
We have isomorphisms of Fr\'echet spaces
$\Cinfa \cong \Cinfc \otimes A$ and $\Cinfap \cong \Cinfcp \otimes A$.
\item
The Fr\'echet algebra $\Cinfa$ is a closed two-sided ideal in $\Cinfap$
and $\Cinfap /\Cinfa \cong A[[z]]$, the Fr\'echet algebra of formal
power series in $z$ with coefficients in $A$, with the topology of
pointwise convergence of the coefficients. If $z(r) = r$ is the identity
function in $\Cinfcp$, then this isomorphism sends
$z^{m} \Cinfap /\Cinfa$ onto $z^{m} A[[z]]$ for each $m\in \N$.
\endroster
\endproclaim
\demo{Proof}
Using the integral form of the remainder for the Taylor
polynomial, for any $C^{\infty}$ function $g\colon [0, 1] \rightarrow
A$ and $r \in [0, 1]$ we obtain
$$ \split \biggl\Vert g(r) - g(0) - &\dots -
 {g^{(k-1)}(0) \over{(k-1)!}}r^{k-1} \biggr\Vert_m \\
& = \biggl\Vert \int_{0}^{r} {g^{(k)}(s) \over{(k-1)!}}(r-s)^{k-1}ds
\biggr\Vert_{m}
\leq {r^k\over {k!}} \sup_{s \in [0, r]} \pa g^{(k)}(s)\pa_{m}.
 \endsplit $$
Let $f\in \Cinfa$, and put $g=f^{(n)}$.  Since the derivatives of
$f^{(n)}$ all vanish at zero, we obtain
$$ \pa r^{-k} \varf^{(n)}(r) \pa_{m}\leq
{1\over{k!}} \sup_{s \in [0, 1]}
\pa \varf^{(n+k)}(s) \pa_{m}.$$
Taking the sup over $r\in [0, 1]$ on the left yields the continuity
of $\pa \quad \pa_{k, n, m}^{\prime}$.
\par
To prove (2), note that the space $\Cinfc$ is a closed subspace of the
nuclear locally convex space $\Cinfcp$ and hence is nuclear
\cite{\Tr, Proposition 50.1}.
\par
For (3), note that the topologies on $\Cinfap$ and $\Cinfa$ are already
given by tensor product seminorms, and the tensor products are
unique by (2).
\par
We prove (4).  Define $\pi \colon \Cinfcp \rightarrow \C[[z]]$
by $\pi(f) = (f(0), f^{(1)}(0), {1\over{2!}}f^{(2)}(0), \dots)$.
Standard differentiation rules show that $\pi$ is a homomorphism,
and  it is trivially continuous.
It is known (see \cite{\Tr, Theorem 38.1}) that for every sequence
$(c_{0}, c_{1}, c_{2}, \dots )$ of complex numbers, there exists
a $C^{\infty}$ complex valued function $f$ on $[0, 1]$ such that
$f^{(n)}(0)= c_{n}$ for all $n\in \N$.  This proves the surjectivity
of $\pi$.
\par
Clearly ${\text{Ker}}(\pi) = \Cinfc$.  Therefore $\pi$ defines a
continuous bijection
${\overline{\pi}} \colon \Cinfcp/\Cinfc \rightarrow \C[[z]]$.
Since both these algebras are Fr\'echet, the open mapping theorem
for Fr\'echet spaces\cite{\Tr, Theorem 17.1}
implies that ${\overline{\pi}}$ is a homeomorphism.
(Note that in \cite{\Tr}, a \lq\lq homomorphism \rq\rq\ is actually
a quotient map; see the beginning of Chapter 17.)
\par
In the notation of the lemma, $\pi(z^{m} \varf) = z^{m}\pi(\varf)$,
and from this one easily sees that ${\overline{\pi}}$
maps $z^{m} \Cinfcp /\Cinfc$ onto $z^{m}\C[[z]]$.
\par
To prove (4), it suffices to show that we can replace $\C$ with $A$.
The map  $ \Cinfap/\Cinfa \rightarrow A[[z]]$
is just the projective tensor product of ${\overline {\pi}}$
with the identity map on $A$.  (Note that $\C [[z]]$ is nuclear
because it is a product of nuclear Fr\'echet spaces
\cite{\Tr, Proposition 50.1}; thus we have
$A[[z]] \cong \C[[z]] \otimes A$.) Since the projective tensor product
of quotient maps of Fr\'echet spaces is a quotient map
\cite{\Tr, Proposition 43.9}, we have the result.
\qed
\enddemo
\subheading{1.1.3 Definition} \cite{\Scha, \S 3}
We say that an action $\alpha$ of
$\R$ on a Fr\'echet algebra $A$ is {\it $m$-tempered} if there exists
a family of submultiplicative seminorms
$\bigl\{\pa \quad \pa_{m}\bigr\}_{m=0}^{\infty}$ giving the topology of
$A$, such that for every $m\in \N$, there is a polynomial $\text{poly}$
such that
 $$ \pa \alpha_{r}(a) \pa_{m} \leq {\text{poly}}(r) \pa a \pa_{m},$$
for $a\in A$ and  $r\in \R$.
(This notion of an $m$-tempered action is a special case of \cite{\Scha,
Definition 3.1.1} with the gauge $\sigma(r) = |r|$ on $\R$.)
For use in \S 2, note that the same definition makes sense for
$\Z$ or $\T$ in place of $\R$.  In the case of $\T$, the polynomial
may be replaced by a constant.
If $\alpha$ is an $m$-tempered action of $\R$ on $A$, then the
Fr\'echet space $\schr \otimes A$ is an ($m$-convex)  Fr\'echet algebra
under the convolution product
$$ (f*g)(r)\,  =\, \int_{\R} \,  f(s)\,  \alpha_{s}(g(r-s))\,   ds, $$
where $f, g \in \schr \otimes A$ and
$ r \in \R$ \cite{\Scha, Theorem 3.1.7}.
We denote this Fr\'echet algebra by $\sona$, and call it the
smooth crossed product.
\proclaim{1.1.4 Lemma} If
$$ 0 \longrightarrow J \longrightarrow A \longrightarrow B
 \longrightarrow 0$$
is a short equivariant exact sequence of Fr\'echet algebras with
$m$-tempered actions $\gamma$, $\alpha$, and  $\beta$ of $\R$, then
$$ 0 \longrightarrow \Cal S(\R, J, \gamma) \longrightarrow \sona
 \longrightarrow \Cal S(\R, B, \beta) \longrightarrow 0$$
is an exact sequence of Fr\'echet algebras.
\endproclaim
\demo{Proof} We can disregard the algebra structure; this then reduces
to the exactness of tensoring with the  nuclear Fr\'echet space
$\schr$ \cite{\Ph, Theorem 2.3(3b)}. \qed \enddemo
\subheading{1.1.5 Definition} Let $\alpha$ be an action of $\R $
on $A$. If for every $a \in A$ the map $r \mapsto \alpha_{r}(a) $
is $C^{\infty}$, we will say that $\alpha$ is a {\it smooth action}.
If $\alpha$ is  smooth and $m$-tempered, we define
the {\it interval action} $\Ialph$ of $\R$ on $\Cinfap$ by
$$\Ialph_{s}(\varf)(r) = \alpha_{rs}(\varf(r))$$
for $s \in \R$, $\varf \in \Cinfap$, and $r\in [0, 1]$.
The {\it cone action} $\Calph$ is defined to be the restriction
of $\Ialph$ to $\Cinfa$.
\proclaim{1.1.6 Lemma} Let $A$, $\alpha$, $\Ialph$, $\Calph$
be as in the previous definition.  Then $\Ialph$ and $\Calph$
are well-defined, smooth, and  $m$-tempered actions of $\R$
on $\Cinfap$ and $\Cinfa$ respectively.
Furthermore, if $z$ is the identity function $z(r)=r$, then for each
$m\in \N$, the subspace $z^{m}\Cinfap$ is a closed $\Ialph$-invariant
ideal in $\Cinfap$.
\endproclaim
\demo{Proof}
Let $\pa \quad \pa_{m}$ be submultiplicative seminorms giving
the topology of $A$.  Let $D$ denote the differential operator
$D\varf (r) = \varf^{(1)}(r)$ for $\varf \in \Cinfap$, and let
${\tilde D}$ denote the differential operator $${\tilde D}\varf (r)
= \lim_{s\rightarrow 0} {\alpha_{s}(\varf(r)) - \varf(r) \over{s}}.
$$
Note that $D$ and $\tilde D$ commute.  The topology on $\Cinfap$ can
then be given by the seminorms
$$ \pa \varf \pa_{m}^{\prime}
= \max_{p+q \leq m} \pa {\tilde D}^{p}{D}^{q}
\varf \pa_{m}
= \max_{p+q \leq m}\biggl( \sup_{r \in [0, 1]}
\pa {\tilde D}^{p}{ D}^{q}\varf(r)\pa_{m}\biggr)$$
The seminorms $\pa \quad \pa_{m}^{\prime}$ are easily seen to be
submultiplicative modulo  constants.  Note that
$$ \pa  {\tilde D}^{p}D^{q} \Ialph_{s}(\varf)\pa_{m}=
\sup_{r\in [0, 1]}
  \pa  {\tilde D}^{p}D^{q} \alpha_{sr}(\varf(r))\pa_{m},$$
where $D$ acts on $r$.  Let $p + q \leq m$.  By the chain rule and since
$\tilde D$ commutes with $\alpha$, this is bounded by
$$ C\sum_{u+t  \leq q} s^{u} \pa \alpha_{sr}
({\tilde D}^{p+u}{D}^{t}\varf(r))\pa_{m}
\leq {\text{poly}}(s)\pa \varf \pa_{m}^{\prime}$$
where we used the $m$-temperedness of the action of $\R$ on $A$.
It follows that {$\pa \Ialph_{s}(\varf)\pa_{m}^{\prime} \leq
{\text{poly}}(s) \pa \varf \pa_{m}^{\prime}$}, so $\Ialph$ is
$m$-tempered.
The smoothness of $\Ialph$ is straightforward using the chain rule.
\par
To prove the statement about $\Cinfa$, it suffices to show that
$\Cinfap$ leaves $\Cinfa$ invariant.  This will follow from
the last statement, since
$\Cinfa = \bigcap_{m=0}^{\infty} z^{m} \Cinfap$.
\par
We therefore prove the last statement. Invariance follows from
$\Ialph (z^{m}\varf) = z^{m}\Ialph (\varf)$.  That $z^{m} \Cinfap$
is an ideal in $\Cinfap$ is obvious.  Closedness follows from
Lemma 1.1.2(4).
\qed
\enddemo
\proclaim{1.1.7 Lemma} Let $\alpha$ be a smooth $m$-tempered action of
$\R$ on a Fr\'echet algebra $A$.  Then
$$ \Cal S(\R, \Cinfa, \Calph) \cong C_{\infty}(\sona)$$  
as Fr\'echet algebras.
\endproclaim
\demo{Proof}
Recall that as a Fr\'echet space, $B_{1}= \Cal S(\R, \Cinfa, \Calph)$ is
isomorphic to $\schr \otimes \Cinfc \otimes A$ with
multiplication       
$$ (f*_{1} g)(s, r) \,=
\,\biggr(\int_{\R} \,f(t) \,\Calph_{t}(g(s-t)) \,dt\,\biggl) (r)\, =
\,\int_{\R} \,f(t, r)\, \alpha_{rt}(g(s-t, r))\,
dt \tag * $$
for $s \in \R$ and  $r\in [0, 1]$.
By commuting $\schr$ and $\Cinfc$ in the tensor product,
we may also identify $B_{2} = C_{\infty}(\sona)$ with $\schr \otimes
\Cinfc \otimes A$, but with multiplication        
$$ (f*_{2} g) (s, r)\, = \,\int_{\R} \,f(t, r)\,
 \alpha_{t}(g(s-t, r))\, dt \tag ** $$
for $s \in \R$ and  $r\in [0, 1]$.
Define $\gamma \colon B_{1} \rightarrow B_{2}$ by $(\gamma f) (s, r)
= r^{-1} f(r^{-1}s, r)$.  We must show this is well defined and
continuous.  By Lemma 1.1.2(1), we may topologize
$\schr \otimes \Cinfc  
\otimes A$ by the seminorms
$$ \pa f \pa_{p, q, l, n, m} = \sup_{s\in\R, r\in [0, 1]}
\w(s)^{p} r^{-q} \pa D^{l} {\tilde D}^{n} f(s,r)\pa_{m},$$
where $\w(s) = 1+ |s|$, and
$D$ and $\tilde {D}$ are derivatives in $s$ and $r$ respectively.
(Note that they commute.)
\par
By the chain rule and the product rule,
$$ \split \pa D^{l}{\tilde D}^{n} &(\gamma f)(s, r) \pa_{m}
= \pa D^{l} {\tilde D}^{n} r^{-1} f(r^{-1}s, r) \pa_{m}\\&
= \pa {\tilde D}^{n} r^{-l-1} (D^{l} f)(r^{-1}s, r) \pa_{m}
\leq C \sum_{i,j,k\leq N} r^{-i}
\pa (D^{j}{\tilde D}^{k} f)(r^{-1}s, r) \pa_{m},
\endsplit $$
for some large $N$ not depending on $f$.  Now if we
multiply a generic summand by $\w(s)^{p} r^{-q}$
and take the sup over $r, s \in \R$, we get something bounded by
$$
\aligned
\sup_{s,r} \w(s)^{p} r^{-q-i}
 \pa (D^{j} {\tilde D}^{k} f) (r^{-1}s, r) \pa_{m} =
\sup_{s,r} \w(sr)^{p} r^{-q-i}&
 \pa (D^{j} {\tilde D}^{k} f) (s, r) \pa_{m} \\
\leq \sup_{s,r} \w(s)^{p} r^{-q-i}
 \pa (D^{j} {\tilde D}^{k} f) (s, r) \pa_{m}
& = \pa f \pa_{p, q+i, j, k, m}, \endaligned
$$
where we replaced $s$ by $sr$ and used $\w(sr) \leq \w(s)$ for
$r \in [0, 1]$. From these estimates,
we see that $\gamma$ is  a well defined map and that
$$ \pa \gamma f \pa_{p,q,k,n,m} \leq C \sum_{i,j,k\leq N}
\pa f \pa_{p, q+i, j, k, m},$$
so $\gamma$ is continuous.  We easily check that $\gamma(f *_{1} g)
= \gamma f *_{2} \gamma g$ using (*) and (**),
so it remains to show that $\gamma$ is a bijection,
and that $\gamma^{-1}$ is continuous.
\par
Define $\tau \colon B_{2} \rightarrow B_{1}$ by $(\tau f)(s, r)
= r f(rs, r)$.  Clearly $\tau$ is an inverse for $\gamma$ as long
as $\tau f$ is in $B_{1}$ for every $f \in B_{2}$.  We have
$$ \pa D^{l}{\tilde D}^{n} (\tau f)(s, r) \pa_{m}
\leq C \sum_{j \leq l, k \leq n}
{\text{poly}}(r) \pa (D^{j}{\tilde D}^{k} f)(sr,r) \pa_{m}.$$
Using $f \in B_{2}$ and $\w(r^{-1}s) \leq r^{-1} \w(s)$,
one easily checks that
$\pa \tau f \pa_{p, q, l, n, m}
 < \infty$ for all $p, q, l, n, m \in \N$.
Hence $\tau f \in B_{1}$ and $\tau$ is an
inverse for $\gamma$.  By the open mapping theorem for Fr\'echet
spaces, $\gamma$ is an isomorphism and we have proved Lemma 1.1.7.
\qed
\enddemo
\subheading{1.1.8 Remark}
Note that the proof above does not work if one uses
$\{ \,\varf \in \Cinfap \,|\, \varf(0)=0 \,\}$ in place of $\Cinfa$.
\vskip\baselineskip
\par
Before we prove our first results on $K$-theory, we prove a lemma on
spectral invariance.  This will be used repeatedly throughout the paper
to show that various dense inclusions induce isomorphisms on $K$-theory.
We say that a subalgebra $A$ of
an algebra $B$ is {\it spectral invariant} in $B$ if
the invertible elements of the unitization $A^{+}$
are precisely the invertible elements of $B^{+}$ which lie in $A^{+}$.
\proclaim{1.1.9 Spectral Invariance Lemma}
\roster
\item
Let $A$ be a dense Fr\'echet subalgebra of a Banach algebra $B$.  If $A$
is spectral invariant in $B$, then the inclusion map
$A\hookrightarrow B$ induces an isomorphism $RK_{*}(A)\cong RK_{*}(B)$.
\item
Let $A$ be a dense Fr\'echet subalgebra of a Fr\'echet algebra
$B$.  Assume that we can write $A$ and $B$ as inverse limits
of Banach algebras $A_{n}$ and $B_{n}$ respectively, where
$A_{n}$ is dense in $B_{n}$ for all $n$,
and the inclusions $A\subseteq A_{n}$, $B\subseteq B_{n}$
are all dense. If each $A_{n}$ is spectral
invariant in $B_{n}$, then the inclusion $A \hookrightarrow B$
induces an isomorphism $RK_{*}(A) \cong RK_{*}(B)$.
\endroster
\endproclaim
\demo{Proof}
For the first item, by Theorem A.2.1 in the appendix of \cite{\Bo},
we know that inclusion induces an isomorphism
$K_{*}(A)\cong K_{*}(B)$.  Since $A^{+}$ and $B^{+}$ have
open groups of invertibles,
$RK_{*}(A) \cong K_{*}(A)$ and $RK_{*}(B) \cong K_{*}(B)$
\cite{\Ph, Theorem 7.7}.  Thus $RK_{*}(A) \cong RK_{*}(B)$.
\par
For the second item,  again by \cite{\Bo, Theorem A.2.1} we have
$K_{*}(A_{n}) \cong K_{*}(B_{n})$ for all $n\in \N$.
Since $A_{n}$ and $B_{n}$ are Banach algebras, their $K_{*}$ is the
same as their $RK_{*}$ \cite{\Ph, Corollary 7.8}.
Applying the Milnor $\varprojlim^{1}$-sequence
\cite{\Ph, Theorem 6.5}, one concludes that $RK_{*}(A)\cong RK_{*}(B)$.
\qed
\enddemo
\par
We let $\schr \otimes A$  denote the Fr\'echet algebra
of Schwartz functions from $\R$ to $A$
with convolution multiplication, for the trivial action
of $\R$ on $A$. (Note that it is isomorphic to the
same set of functions with pointwise multiplication via the
Fourier transform.)  We find $\schr \otimes A$ a more convenient
suspension than, say, $C_{0}(\R, A)$.  Part (3) of the
next lemma shows that it behaves the same way for $K$-theory.
\proclaim{1.1.10 Lemma} Let $A$ be a Fr\'echet algebra.  Then:
\roster
\item Evaluation at $r\in [0, 1]$ induces an isomorphism
$$ ({\text{ev}}_{r})_{*} \colon RK_{*}(\Cinfap) \rightarrow RK_{*}(A), $$
which does not depend on $r$.
\item $RK_{*}(\Cinfa)=0$.
\item There is a natural isomorphism
$RK_{*}(A)\cong  RK_{*+1}(\schr \otimes A)$.
\endroster
\endproclaim
\demo{Proof}
(1) The maps ${\text{ev}}_{r}$ are all homotopic, so the maps
$({\text{ev}}_{r})_{*}$ are all equal.  To see that they are
isomorphisms, it suffices to show that one of them, say
${\text{ev}}_{0}$, is a homotopy equivalence.  The homotopy
inverse $\varphi \colon A \rightarrow \Cinfap$ is given by
$\varphi(a)(r) = a$.
The homotopy from $\varphi \circ {\text{ev}}_{0}$
to ${\text{id}}_{\Cinfap}$is given by
$\psi_{t}(\varf)(r)= \varf(tr)$ for $t, r \in [0, 1]$.
One easily verifies that everything is in fact smooth and continuous
in the $C^{\infty}$ topology.
\par
(2) The identity map of $\Cinfa$ is homotopic to the zero map via
$\psi_{t}(\varf)(r) = \varf(tr)$.
\par
(3) Let $\bigl\{\pa \quad \pa_{m}\bigr\}_{m=0}^{\infty}$
be an increasing sequence of submultiplicative seminorms
which define the topology of $A$, and let $A_{m}$ denote the
Banach algebra obtained by completing $A/{\text{Ker}}(\pa \quad
\pa_{m})$.  Set
$$\split B_{m} = \biggl\{\, \varf \in & C^{m}(\R, A_{m})\,\biggl|\\
& \pa f \pa_{m}^{\prime} = \max_{0 \leq k \leq m}
\biggl(\sup_{s\in \R} (1+|s|)^{m} |f^{(k)}(m)|
\biggr)
<\infty\,\biggr\}.\endsplit$$
Then $B_{m}$ with pointwise multiplication is a Banach algebra
in a norm equivalent to $\pa \quad \pa^{\prime}_{m}$,
and $B_{m}$ is a dense subalgebra of the Banach algebra $C_{m}=
C_{0}(\R, A_{m})$ of continuous functions from $\R$ to
$A_{m}$ which vanish at infinity.  It is easily checked that
$B_{m}$ is spectrally invariant in $C_{m}$\cite{\Schb, (2.3)-(2.4)}.
Also,  $\varprojlim C_{m} \cong C_{0}(\R, A)$, which is the usual
suspension
of $A$, and  $\varprojlim B_{m} \cong \schr \otimes A$, where $\schr$
is taken with pointwise multiplication.  The Fourier transform shows
that this algebra is isomorphic to the one obtained using convolution on
$\schr$.  Therefore, using Lemma 1.1.9 (2) in the first step and Bott
periodicity \cite{\Ph, Theorem 5.5} in the second step, we obtain
$$ RK_{*+1}(\schr \otimes A) \cong RK_{*+1} (C_{0}(\R, A))
 \cong RK_{*}(A).$$
\qed \enddemo
\par
Item (3) in the previous lemma can also be proved using homotopy
invariance results of A.M. Davie
\cite{\Da}, instead of using spectral invariance.
\par
We let $\Kinf$  denote the Fr\'echet algebra of
rapidly vanishing complex valued matrices on $\Z \times \Z$.
This is the smooth version of the compact operators used in
Section 2 of \cite{\Ph} and \cite{\Schb, \S 5}.
\proclaim{1.1.11 Lemma} Let $\alpha$ be a smooth $m$-tempered
action of $\R$ on the Fr\'echet algebra $A$. Let $\beta$ be the
action on the quotient $\Cinfap/\Cinfa\cong A[[z]]$ determined by
$\Ialph$.  Then
${\text{ev}}_{0} \colon \Cinfap \rightarrow A$ determines an
isomorphism
$$ RK_{*}(\Cal S(\R, \Cinfap/\Cinfa, \beta))
\cong RK_{*}(\schr \otimes A).$$
\endproclaim
\demo{Proof}
Note that ${\text{ev}}_{0}$
is equivariant for the action $\Ialph$ on $\Cinfap$
and the trivial action $\tau$ on $A$.  Therefore we obtain
$${\overline{{\text{ev}}}}_{0} \colon \Cinfap/\Cinfa \rightarrow A$$
and
$${\overline{\epsilon}}_{0} \colon
\Cal S(\R, \Cinfap/\Cinfa, \beta) \rightarrow \Cal S(\R, A, \tau)
\cong \schr \otimes A.$$
We have to show that $({\overline{\epsilon}}_{0})_{*}$ is an
isomorphism.   By Lemma 1.1.6 and Lemma 1.1.2(4), we have
$\Cinfap/\Cinfa \cong A[[z]]$, and the induced action $\gamma$
of $\R$ on $A[[z]]$ satisfies $\gamma_{s}(z^{m}A[[z]])
\subseteq z^{m} A[[z]]$ for all $m\in \N$.  It is easily verified
that if $a= a_{0} + a_{1} z + \dots \in A[[z]]$, then
$\gamma(a) = b_{0} + b_{1}z + \dots $ satisfies
$b_{0}= a_{0}$.  Thus $a\mapsto a_{0}$ is an equivariant
map from $(A[[z]], \gamma)$ to $A$ with the trivial action.
\par
Let $B = \Cal S(\R, A[[z]], \gamma)$, and let $B_{m}$
be the closed ideal $\Cal S(\R, z^{m}A[[z]], \gamma)$.
Then $B=B_{0}$, and $B/B_{1} \cong SA$.  We have to prove that the
quotient map $B \rightarrow B/B_{1}$ is an isomorphism on $K$-theory.
Considering the long exact sequence in $K$-theory
\cite{\Ph, Theorem 6.1}, it suffices to show that $RK_{*}(B_{1})=0$.
\par
We first prove $RK_{1}(B_{1})=0$.  This is equivalent to showing that
the invertible group ${\text{inv}}((\Kinf \otimes B_{1})^{+})$
is connected.  As a topological vector space (TVS), we have
$zA[[z]] \cong \prod_{m=1}^{\infty} z^{m}A[[z]]/
z^{m+1}A[[z]]$.  Further, as TVS's, we have
$\Cal S(\R, D, \delta) \cong \schr \otimes D$ for any $D$ and
$\delta$.  Hence $B_{1} \cong \prod_{m=1}^{\infty} B_{m}/B_{m+1}$
as TVS's. It follows that $\Kinf \otimes
B_{1} \cong \prod_{m=1}^{\infty} (\Kinf \otimes B_{m})/(\Kinf \otimes
B_{m+1})$
as TVS's.
\par
For $m, n\in \N$, we have
$(z^{m} A[[z]])(z^{n}A[[z]]) \subseteq z^{m+n}A[[z]]$.
The fact that $z^{m+n}A[[z]]$ is closed now implies that
$B_{m} B_{n} \subseteq B_{m+n}$, so
$(\Kinf \otimes B_{m})(\Kinf \otimes B_{n}) \subseteq
\Kinf \otimes B_{m+n}$.
\par
We will now show that if $\eta \in \C - \{0\}$ and $x \in \Kinf \otimes
B_{1}$ then $\eta -x$ is invertible.  This will clearly imply that
${\text{inv}} ((\Kinf \otimes B_{1})^{+})$ is connected.  Without loss
of generality, let $\eta = 1$.  Note that $x^{n} \rightarrow 0$ in
$\Kinf \otimes B_{1}$, because the image of $x^{n}$ in
$ \prod_{m=1}^{\infty} (\Kinf \otimes B_{m})/(\Kinf \otimes B_{m+1})$
is zero in the first $n-1$ factors.  Furthermore, the sum
$\sum_{n=1}^{\infty} x^{n}$  converges in $\Kinf \otimes B_{1}$,
since $\sum_{n=N}^{N+p} x^{n}$ is zero in the first $N-1$ factors of
$\prod_{m=1}^{\infty} (\Kinf \otimes B_{m})/(\Kinf \otimes B_{m+1})$,
and the topology is given by pointwise convergence in each factor.
Therefore $1-x$ is invertible with inverse $\sum_{n=1}^{\infty}x^{n}$.
\par
It remains to show that $RK_{0}(B_{1})=0$.  If we replace $A$ by
$C^{\infty}(\T)\otimes B$, it is easy to see (using the nuclearity
of $C^{\infty}(\T)$) that this replaces $B_{1}$ by $C^{\infty}
(\T) \otimes B_{1}$.   Therefore we know from the above that $RK_{1}
(C^{\infty}(\T)\otimes B_{1})=0$.
By the estimates \cite{\Schb, (2.3)-(2.4)}, the inclusion
$C^{\infty}(\T) \otimes B_{1} \hookrightarrow C(\T, B_{1})$
is an isomorphism on $K$-theory.  Let $SB_{1} = \{ f \in C(\T, B_{1}) |
f(1)=0\}$.  Then from the split exact sequence
$$ 0 \rightarrow SB_{1} \rightarrow C(\T, B_{1}) \rightarrow B_{1}
\rightarrow 0, $$
and Bott periodicity \cite{\Ph, Theorem 5.5}, we have
$ RK_{1}(C(\T, B_{1})) = RK_{0}(B_{1}) \oplus RK_{1}(B_{1})$.
Hence $RK_{0}(B_{1})=0$.
\qed
\enddemo
\proclaim{1.1.12 Corollary} The map
$$ \epsilon_{0} \colon \Cal S(\R, \Cinfap, \Ialph) \rightarrow
\schr \otimes A,$$
induced on the crossed products by ${\text{ev}}_{0} \colon \Cinfap
\rightarrow A$, is an isomorphism on $K$-theory.
\endproclaim
\demo{Proof}
Using the notation of the previous lemma, we can factor $\epsilon_{0}$
as
$$ \Cal S(\R, \Cinfap, \Ialph)\rightarrow
\Cal S(\R, \Cinfap/\Cinfa, \beta) \rightarrow \schr \otimes A.$$
The kernel of the first map is $\Cal S (\R, \Cinfa, \Calph)$
(by Lemma 1.1.4), so Lemma 1.1.7 and Lemma 1.1.10 (2)
imply that its $K$-theory is $0$.  The long exact sequence for
$K$-theory \cite{\Ph, Theorem 6.1} therefore shows that the first map
is an isomorphism on $K$-theory.  The second map is an isomorphism on
$K$-theory by Lemma 1.1.11.
\qed
\enddemo
\vskip\baselineskip
\heading \S 1.2 The Thom Map $\theta$ is an Isomorphism \endheading
\par
We let $\beta$ (or, when necessary, $\beta_{A}$) denote
the Bott periodicity isomorphism $RK_{*}(A)
\rightarrow RK_{*+1}(\schr \otimes A)$, as obtained in Lemma
1.1.10(3).
\par
Throughout this section, we shall use $I(\alpha)$ to denote
$\Ialph$.
\subheading{1.2.1 Definition} Let $\alpha$ be a smooth
$m$-tempered action of $\R$ on a Fr\'echet algebra $A$.
We define $F(A, \alpha) =\Cal S(\R, \Cinfap, I(\alpha))$.
(This corresponds to the algebra of sections of the continuous
field of \cite{\ENNb}.) We let
$$\epsilon_{0} \colon F(A, \alpha) \rightarrow \schr \otimes A$$
be as in Corollary 1.1.12, and we let $$\epsilon_{1}
\colon F(A, \alpha) \rightarrow \sona$$ be the map on the crossed
products determined by ${\text{ev}}_{1} \colon \Cinfap \rightarrow A$.
We then define
$$ \theta^{(0)}\colon RK_{*}(\schr \otimes A)
\rightarrow RK_{*}(\sona)$$
by $\theta^{(0)} = (\epsilon_{1})_{*}\circ (\epsilon_{0})_{*}^{-1}$.
(We sometimes write $\theta_{A}^{(0)}$ or $\theta^{(0)}_{A, \alpha}$.)
We further define the Thom map
$$ \theta\colon RK_{*}(A) \rightarrow RK_{*+1}(\sona)$$
by $\theta = \theta^{(0)} \circ \beta$
(or $\theta_{A} = \theta_{A}^{(0)} \circ \beta_{A}$).
\proclaim{1.2.2 Lemma} The  maps $\theta^{(0)}$ and $\theta$
are well defined natural group homomorphisms.
Moreover, they both commute with Bott periodicity.  That is, the diagram
$$\CD
RK_{*}(\schr \otimes A) @>{\theta_{A}^{(0)}}>>
RK_{*}(\sona) \\
@V{\beta}VV    @VV\beta V\\
RK_{*+1}(\schr \otimes \schr \otimes A) @>{\theta_{\schr
 \otimes A}^{(0)}}>>
RK_{*+1}(\Cal S(\R, \schr \otimes A, \alpha)) \endCD
$$
commutes, and similarly for $\theta$ in place of $\theta^{(0)}$.
\endproclaim
\demo{Proof}
Since $\theta$ is essentially $\theta^{(0)} \circ \beta$, we
need only prove this for $\theta^{(0)}$.  Since $\theta^{(0)}
= (\epsilon_{1})_{*} \circ (\epsilon_{0})_{*}^{-1}$, this follows from
the naturality of Bott periodicity and the isomorphism
$F(\schr \otimes A, {\text{id}}\otimes
\alpha ) \cong \schr \otimes F(A, \alpha)$.\qed \enddemo
\proclaim{1.2.3 Lemma} Let $\alpha$ be an $m$-tempered action
of $\R$ on a Fr\'echet algebra $A$.  Then the dual action
$\widehat{ \alpha}$ of $\R$ on $\sona$ is smooth and $m$-tempered.
\endproclaim
\demo{Proof}
Recall that ${\widehat{\alpha}}_{s}(\varf)(r) = e^{2\pi i rs} \varf(r)$,
where $s, r \in \R$ and $\varf \in \sona$.  Since differentiating
in $s$ only brings down powers of $2\pi ir$ and $\sona$ consists of
Schwartz functions, we see that $\widehat{\alpha}$ is a smooth action.
\par
Let $\bigl\{\pa \quad \pa_{m}\bigr\}$ be a family of submultiplicative
seminorms giving the topology on  the convolution algebra of
$L^{1}$-rapidly vanishing functions $L_{1}^{|\,\,|}(\R, A, \alpha)$ from
$\R$ to $A$ \cite{\Scha, (2.1.3), Thm 3.1.7}.
(The superscript $|\,\, |$ means that the functions are integrable
against any power of the absolute value function $|\,\, | $.)
Let $B= L_{1}^{|\,\,|}(\R, A, \alpha)$.
We show that $\pa \quad \pa_{m}$ are in fact equivalent to
a family of submultiplicative seminorms on which $\alphh$ acts
isometrically on each seminorm.  It is easy to see that the seminorms
$$\pa \varf\pa_{m}^{\prime}=
 \int_{\R} (1+|r|)^{m} \pa \varf(r) \pa_{m} dr, $$
where $\varf \in B$, are isometric for the action
$\alphh$ (since $\pa e^{2\pi i rs} \varf(r) \pa_{m} =
 \pa \varf (r) \pa_{m}$).
These norms $\pa\quad\pa_{m}^{\prime}$ topologize
$B$, and so are equivalent to
the $\pa \quad \pa_{m}$.  Hence for each $\varf \in B$,
the function $s\mapsto \pa\alphh_{s}(\varf)\pa_{m}$ is bounded
by  $C\pa \varf \pa^{\prime}_{k}$ for some sufficiently
large $k\in\N$, $C>0$ chosen independently of $\varf$.
Then the functions $f\mapsto
\sup_{s\in \R}\pa \alphh_s(\varf) \pa_{m}	$
are easily seen to be a well-defined family of submultiplicative,
$\alphh$-isometric, seminorms topologizing $B$.  Replace our
original seminorms $\pa \quad \pa_{m}$ with these ones.
\par
The seminorms
$$ \pa \varf \pa_{m}^{\prime} =
\max_{k \leq m} \pa D^{k}\varf \pa_{m}$$
topologize $\sona$ and are submultiplicative since
$$  \pa \varf * g \pa_{m}^{\prime} =
\max_{k \leq m} \pa \varf * D^{k}g \pa_{m} \leq
\pa \varf \pa_{m} \max_{k \leq m}\pa D^{k}g\pa_{m}
\leq \pa \varf \pa_{m}^{\prime} \pa
g \pa_{m}^{\prime},
$$
for $f,g \in \sona$.  By the product rule,
$$ \pa D^{k}{\widehat{\alpha}}_{s}(\varf) \pa_{m} \leq
{\text{poly}}(s)\biggl( \max_{i \leq k} \pa\alphh_{s} D^{i} \varf\pa_{m}
\biggr)
= {\text{poly}}(s) \biggl(\max_{i \leq k} \pa D^{i} \varf\pa_{m}
\biggr),$$
so $\pa \alphh_{s}(\varf)\pa_{m}^{\prime}
\leq {\text{poly}}(s)\pa \varf \pa_{m}^{\prime}$
and $\alphh$ is $m$-tempered (see Definition 1.1.3).
\qed
\enddemo
\subheading{1.2.4 Remark} Note that the above proof shows that
it is possible to have an $m$-tempered action of $\R$
on a Fr\'echet algebra $A$ for which there is no choice of
seminorms (submultiplicative or otherwise) for which $\alpha$
acts isometrically.  For example, let $A=\schr$ with
pointwise (convolution) multiplication, and let $\R$ act by
translation (dual of translation).  Then for any family
$\bigl\{ \pa \quad \pa_{m} \bigr\}$ of seminorms topologizing
$\schr$, and for any nonzero $\varf \in \schr$,
there must be infinitely many $m$ such that
$\pa \alpha_{s}(\varf)\pa_{m}$ is not  bounded in $s$.
\vskip\baselineskip
\par
We remark that the $m$-temperedness of $\widehat \alpha$ is not needed
for the $m$-convexity of the crossed product
$\Cal S(\R, \sona, {\widehat \alpha})$
by the following lemma ($m$-temperedness is not required in its proof).
However Lemma 1.2.3 will be needed, for example,  for the $m$-convexity
of the algebra $F(F(A, \alpha),{\widehat {I(\alpha)}})$ used below.
\proclaim{1.2.5 Lemma} Let $A$, $\alpha$ and ${\widehat{\alpha}}$
be as in the previous lemma.  Then there is a natural
isomorphism (Takesaki-Takai duality)
$$d_{\alpha} \colon \Cal S( \R, \sona, \alphh ) \cong \Kinf \otimes A.$$
\endproclaim
\demo{Proof} This is proved without the $m$-convexity conditions
in Lemma 2.8 of \cite{\ENNa}, but with a different definition of
$\Kinf$ than ours.  Define
$ \Kinf_{\R} =  \Cal S(\R^{2}) $ with matrix multiplication
$(f*g)(r, t) = \int_{\R} f(r, s) g(s, t) ds$.
Let $\bigl\{\xi_{n}\bigr\}_{n \in \N}$ be the orthonormal basis of
$L_{2}(\R)$
consisting of the Hermite functions. For $f \in \Cal S(\R)$ define
constants $c_{n}$ by the formula  $f= \sum_{n \in \N} c_{n}\xi_{n}$,
and define a map $\schr \rightarrow \Cal S(\N)$ by
$f \mapsto \{ c_{n} \}_{n \in \N}$.  This gives an inner product
preserving isomorphism of $\schr$ with $\Cal S(\N)$
\cite{\RS, Theorem V.13}.
The tensor product of this map with itself gives an isomorphism
$\schr \otimes \schr \rightarrow \Cal S(\N) \otimes \Cal S(\N)$,
which is a Fr\'echet algebra isomorphism of $\Kinf_{\R}$ with $\Kinf$.
\qed
\enddemo
\proclaim{1.2.6 Lemma} Let $\alpha$  be a smooth, $m$-tempered
action of $\R$ on a Fr\'echet algebra $A$. Let $\tau$ be the trivial
action.  Then the following diagram commutes:
$$\minCDarrowwidth{1pt}
\CD
RK_{*}(\schr \otimes \schr \otimes A)
  @>{{\text{id}}\otimes \theta^{(0)}_{A,\alpha}}>>
RK_{*}(\schr\otimes \sona)  @>{\theta^{(0)}_{\sona, {\hat \alpha}}}>>
RK_{*}(\Cal S(\R, \sona, \alphh)) \\
@V{ {\text{id}}\otimes \theta^{(0)}_{A, \tau}}VV
@.
@VV{\scriptstyle (d_{\alpha})_{*}}V \\
RK_{*}(\schr \otimes\Cal S(\R, A, \tau))
@>{\theta^{(0)}_{S(\R, A, \tau), {\hat \tau}}}>>
RK_{*}(\Cal S(\R,  \Cal S(\R, A, \tau), {\hat{\tau}}))
@>{(d_{\tau})_{*}}>>
RK_{*}(\Kinf \otimes A) \endCD $$
\endproclaim
\demo{Proof}
We make the obvious identification of $\Cal S(\R, A,\tau)$
with $\schr \otimes A$. (Since $I_{\infty}(\tau)$ is the trivial action,
$F(A, \tau) \cong \schr \otimes \Cinfap$, and the left vertical
map is induced by the identity.)
Consider the following diagram of Fr\'echet algebras.
$$\CD
\schr \otimes \schr \otimes A @<{{\text{id}}\otimes \epsilon_{0}}<<
\schr  \otimes F(A, \alpha) @>{{\text{id}} \otimes \epsilon_{1}}>>
\schr \otimes \sona \\
@A{\epsilon_{0}}AA
@A{\epsilon_{0}}AA
@A{\epsilon_{0}}AA
\\
F(\schr\otimes A, {\widehat{\tau}}) @<{F(\epsilon_{0})}<<
F(F(A, \alpha), {\widehat{I(\alpha)}}) @>{F(\epsilon_{1})}>>
F(\sona, \alphh)\\
@V{\epsilon_{1}}VV
@V{\epsilon_{1}}VV
@V{\epsilon_{1}}VV
\\
\Cal S(\R, \schr \otimes A, {\widehat{\tau}})
@<{{\overline{ \epsilon}}_{0}}<<
\Cal S(\R, F(A, \alpha),
{\widehat{I(\alpha)}} )
@>{{\overline{ \epsilon}}_{1}}>>
\Cal S(\R, \sona, \alphh) \\
@V{ d_{\tau}}VV
@V{d_{I(\alpha)}}VV
@V{d_{\alpha}}VV
\\
\Kinf\otimes A @<{{\text{id}} \otimes {\text{ev}}_{0}}<<
\Kinf \otimes \Cinfap @>{{\text{id}}\otimes {\text{ev}}_{1}}>>
\Kinf \otimes A
 \endCD \tag * $$
Here ${\overline{\epsilon}}_{0}$ and ${\overline{\epsilon}}_{1}$
are the maps induced on the crossed products by $\epsilon_{0}$ and
$\epsilon_{1}$.
\par
Naturality of the construction $F$ and of the maps $\epsilon_{0}$
and $\epsilon_{1}$ implies that the top four squares  commute,
and naturality of Takesaki-Takai duality gives commutativity
for the bottom two squares.  All the upward arrows along the top
(labelled $\epsilon_{0}$) are isomorphisms on $K$-theory by
Corollary 1.1.12. We claim that the same is true for the left arrows
on the left side.  This is certainly the case for
${\text{id}} \otimes \epsilon_{0}$,
since it is the suspension of an isomorphism on $K$-theory.  Therefore
$F(\epsilon_{0})_{*}$ is an isomorphism, because the other three
maps in the top left square are isomorphisms on $K$-theory.
Furthermore,
${\text{ev}}_{0}$
is an isomorphism on $K$-theory by Lemma 1.1.10(1), so its
stabilization ${\text{id}}
\otimes {\text{ev}}_{0}$ is too. Since $d_{\tau}$ and
$d_{I(\alpha)}$ are actually algebra isomorphisms, it follows
that $({\overline{\epsilon}}_{0})_{*}$ is an isomorphism.  This
proves the claim.
\par
We now apply $RK_{*}$ to the entire diagram $(*)$.  The $K$-theory
diagram still commutes when all the top vertical arrows
and the left horizontal arrows are inverted.  Since
$\theta^{(0)} = (\epsilon_{1})_{*} \circ (\epsilon_{0})_{*}^{-1}$,
we obtain the commutative diagram
$$\CD
RK_{*}(\schr \otimes \schr \otimes A) @>{{\text{id}}\otimes
\theta^{(0)}_{A,\alpha}}>>
RK_{*}(\schr \otimes  \sona) \\
@V{\theta^{(0)}_{\schr \otimes A,{\hat{\tau}}}}VV
@VV{\theta^{(0)}_{\sona, {\hat \alpha}}}V
\\
RK_{*}(\Cal S(\R,  \schr \otimes A, {\widehat{\tau}})) @.
RK_{*}(\Cal S(\R,\sona, \alphh))\\
@V{{(d_{\tau})_{*}}}VV
@VV{{(d_{\alpha})_{*}}}V \\
RK_{*}(\Kinf \otimes A)
@>{({\text{id}}\otimes {\text{ev}}_{1})_{*} \circ
({\text{id}} \otimes {\text{ev}}_{0})_{*}^{-1}}>>
RK_{*}(\Kinf \otimes A)\endCD. $$
The bottom horizontal map is the identity, so this is the diagram
in the statement of the lemma.
\qed
\enddemo
\proclaim{1.2.7 Theorem} Let $\alpha$ be a smooth $m$-tempered action
of $\R$ on a Fr\'echet algebra $A$.  Then the Thom map $\theta \colon
RK_{*}(A) \rightarrow RK_{*+1}(\sona)$ is an isomorphism.
\endproclaim
\demo{Proof} Let $\tau$ denote the trivial action of $\R$ on $A$.
Consider the following diagram:
$$\minCDarrowwidth{1pt}
\CD
RK_{*}(A)    @.
{\theta_{A}\qquad \qquad \qquad}
@.\\
@V\beta_{A}VV
@. @. \\
RK_{*+1}(\schr \otimes A)  @>{\theta^{(0)}}_{A, \alpha}>>
RK_{*+1}(\sona)
@.
{\theta_{\sona, {\hat {\alpha}}}\qquad\qquad\qquad}
\\
@V{\beta_{\schr \otimes A}}VV
@VV{\beta_{\sona}}V
@.
\\
RK_{*}(\schr \otimes \schr \otimes A)
  @>{{\text{id}}\otimes\theta^{(0)}_{A,\alpha}}>>
RK_{*}(\schr\otimes \sona)  @>{\theta^{(0)}_{\sona, {\hat \alpha}}}>>
RK_{*}(\Cal S(\R, \sona, \alphh)) \\
@V{ {\text{id}}\otimes\theta^{(0)}_{A, \tau}}VV
@.
@VV{\scriptstyle (d_{\alpha})_{*}}V \\
RK_{*}(\schr \otimes\Cal S(\R, A, \tau))
@>{\theta^{(0)}_{S(\R, A, \tau), {\hat \tau}}}>>
RK_{*}(\Cal S(\R,  \Cal S(\R, A, \tau), {\widehat{\tau}}))
@>{(d_{\tau})_{*}}>>
RK_{*}(\Kinf \otimes A) \endCD $$
The two triangles commute by definition, the square on the
middle left commutes by Lemma 1.2.2, and the rectangle at the bottom
commutes by Lemma 1.2.6.  Thus the diagram commutes.
\par
We now claim that the map from $RK_{*}(A)$ to $RK_{*}(\Kinf \otimes A)$
is the usual stabilization isomorphism.  We will prove this by
going along the left side and the bottom, and reducing the
problem to Theorem 4.5 of \cite{\ENNb}.  Using Bott periodicity
it suffices to prove this for $RK_{0}$.  Let $\eta \in RK_{0}(A)$,
and choose a representative idempotent $p \in B = M_{2}((\Kinf
\otimes A)^{+})$.  Note that
$p -\pmatrix 1 & 0 \\ 0 & 0 \endpmatrix \in M_{2}(\Kinf \otimes A)$.
Define $\varphi \colon \C \rightarrow M_{2}((\Kinf
\otimes A)^{+})$ by $\varphi(\lambda) = \lambda p$.
Then we have the following commutative diagram.
$$\minCDarrowwidth{15pt}
\CD
RK_{0}(A) @>{\beta\circ \beta}>> RK_{0}(\schr \otimes \schr \otimes A)
@>{\theta^{(0)}}>>
RK_{0}(\Cal S( \R, \schr \otimes A, {\widehat{\tau}}))
@>{(d_{\tau})_{*}}>>  RK_{0}(\Kinf \otimes A) \\
@V{\psi_{*}}VV
@V{({\text{id}}\otimes {\text{id}}\otimes \psi)_{*}}VV
@VVV
@V{({\text{id}} \otimes \psi)_{*}}VV\\
RK_{0}(B) @>{\beta \circ \beta}>>
RK_{0}(\schr \otimes \schr \otimes B) @>{\theta^{(0)}}>>
RK_{0}(\Cal S( \R, \schr \otimes B, {\widehat{\tau}}))
@>{(d_{\tau})_{*}}>>
RK_{0}(\Kinf \otimes B)
\\
@A{\varphi_{*}}AA    @A{({\text{id}}\otimes
{\text{id}}\otimes\varphi)_{*}}AA
@AAA
@A{({\text{id}} \otimes \varphi)_{*}}AA \\
RK_{0}(\C) @>{\beta \circ \beta}>>
RK_{0}(\schr \otimes \schr \otimes \C) @>{\theta^{(0)}}>>
RK_{0}(\Cal S(\R, \schr, {\widehat{\tau}}))
@>{(d_{\tau})_{*}}>>  RK_{0}(\Kinf)\\
@VVV @VVV @VVV @VVV
\\
K_{0}(\C)
@>{\beta \circ \beta}>>
K_{0}(C_{0}(\R^{2}))
@>{\theta^{*}}>>
K_{0}(C^{*}(\R, C_{0}(\R), {\widehat{\tau}})) @>>>
K_{0}(\Cal K)\endCD
$$
Note that the composition of the horizontal maps on the top
is precisely the map from $RK_{*}(A)$ to $RK_{*}(\Kinf \otimes A)$
along the left and bottom of the previous diagram.
The top vertical arrows are induced by the map $RK_{0}(A) \cong
RK_{0} (M_{2}(\Kinf \otimes A)) \rightarrow RK_{0}(B)$ coming from
$\psi \colon M_{2}(\Kinf \otimes A) \hookrightarrow B$.
The bottom vertical arrows are induced by inclusions
in the corresponding C*-algebras, and the groups along the bottom row
are the conventional $K$-theory of C*-algebras.  The map $\theta^{*}$
in the bottom row is constructed in the same manner as $\theta^{(0)}$.
That is, let $D=C_{0}(\R^{2})$, let
${\text{ev}}_{0}, {\text{ev}}_{1} \colon
C([0, 1]) \otimes D \rightarrow D$ be the evaluation maps, and let
$\epsilon_{0}$ and $\epsilon_{1}$ be the crossed products of
${\text{ev}}_{0}$
and ${\text{ev}}_{1}$ by $\R$; then $\theta^{*} = (\epsilon_{1})_{*}
\circ (\epsilon_{0})^{-1}_{*}$.
It is much easier
to prove that $(\epsilon_{0})_{*}$ is invertible in the C*-algebra
case, since one can use $C_{0}((0, 1]) \otimes D$ in place of
$C_{\infty}D$.  
It is also trivial to check that $\theta^{*}$ is the same map as in
Theorem 4.5 of \cite{\ENNb}.
\par
It follows from this same theorem that the map $K_{0}(\C)
\rightarrow K_{0}(\Cal K)$ along the bottom row is the
usual stabilization map.  Now the maps $RK_{0}(\C) \rightarrow
K_{0}(\C)$ and $RK_{0}(\Kinf ) \rightarrow K_{0}(\Cal K)$
are isomorphisms, the range of the map $\varphi_{*}$ contains
$[p] = \varphi_{*}([1])$, the maps $\psi_{*}$ and
$({\text{id}} \otimes \psi)_{*}$ are injective since they come from
split exact sequences, and $[p] = \psi_{*}([p])$ is in the
range of $\psi_{*}$.  Therefore the map
$RK_{0}(A)\rightarrow RK_{0}(\Kinf \otimes A)$ across the top row
sends $\eta = [p]$ to its image under the usual stabilization
map.  Since $\eta \in RK_{0}(A)$ was arbitrary, this proves the claim.
\par
It follows that $\theta_{\sona} \circ \theta_{A}$ is an isomorphism.
Applying this to the action $\alphh$ of $\R$ on $\sona$ now shows that
$\theta_{\sona}$, and hence also $\theta_{A}$, is an isomorphism.
\qed\enddemo
\heading \S 1.3 Replacing $A^{\infty}$ with $A$ and $\schr$ with
 $L_{1}(\R)$ \endheading
\par
We show that the smooth Thom isomorphism still holds if the action of
$\alpha$ on $A$ is strongly continuous but not smooth.  We also
show that if $\alpha$ leaves each seminorm on $A$ invariant, then the
Schwartz algebra $\sona$ may be replaced by $L_{1}(\R, A, \alpha)$.
\proclaim{1.3.1 Theorem} Let $A$ be a Fr\'echet algebra with
$m$-tempered and strongly continuous action $\alpha$ of $\R$.
Let $A^{\infty}$ denote the Fr\'echet algebra of $C^{\infty}$-vectors.
Then the inclusion map
$$ \Cal S(\R, A^{\infty}, \alpha) \hookrightarrow \sona $$
induces an isomorphism on $K$-theory. \endproclaim
\demo{Proof}
If $f \in \sona$, define $\alpha_{s}(f) (r) = \alpha_{s}(f(r))$.
Then $\alpha$  defines an action of $\R$ on $\sona$ by
{\it algebra} automorphisms.  Viewing $\sona$ as
the Fr\'echet space $\schr \otimes A$, note that $\alpha$ is just the
tensor product of the identity map on $\schr$ with the original
action of $\alpha$ on $A$.  The strong continuity of
$\alpha$ on $\sona$ follows easily from the strong continuity
of $\alpha$ on $A$.  The set of $C^{\infty}$-vectors for the
action of $\alpha$ on $\sona$ is easily seen to be
$\Cal S(\R, A^{\infty}, \alpha)$.
(For example, use \cite{\Scha, Theorem A.8}
and the nuclearity of $\schr$.)
\par
We show that $\alpha$ is an $m$-tempered action of $\R$ on $\sona$.
To do this, we will show that the condition for an $m$-tempered action
introduced in Theorem 3.1.18 of \cite{\Scha} is satisfied.
Since the action of $\alpha$ on $A$ is $m$-tempered by assumption,
we let $\pa \quad \pa_{k}$
be an increasing family of submultiplicative seminorms topologizing
$A$, satisfying the conditions of Definition 1.1.3.
For fixed $k$, let $C>0$ and $d \in \N$ be such that
$$ \pa \alpha_{s}(a) \pa_{k} \leq C\w(s)^{d}
\pa a \pa_{k}, \qquad s\in \R,\quad a \in A, \tag * $$
where $\w(s) = 1+|s|$.
Define seminorms (not necessarily submultiplicative) for $\sona$ by
$$ \pa f \pa_{m, l, k} =
\int_{\R} \w(r)^{m} \pa f^{(l)}(r) \pa_{k} dr, \qquad f \in \sona.$$
\par
Let $s_{1}, \dots s_{n} \in \R$ and $f_{1}, \dots f_{n} \in \sona$.
In the following computation, we let
$r_{n}= r-r_{1}-\dots -r_{n-1}$, $\eta_{1}=0$, and
$\eta_{k}= r_{1}+\dots +r_{k-1}$. Then
$$ \aligned &
\pa  \alpha_{s_{1}}(f_{1}) \dots \alpha_{s_{n}}(f_{n}) \pa_{m, l, k}
 \\
&
=\int_{\R}  \w(r)^{m} \biggl\Vert\int_{\R}
\dots \int_{\R} \alpha_{s_{1}+\eta_{1}}(f_{1}(r_{1}))
\\
&\qquad \qquad \qquad \qquad \qquad
\dots
\alpha_{s_{n-1} + \eta_{n-1}} (f_{n-1}(r_{n-1}))
\alpha_{s_{n}+\eta_{n}}(f^{(l)}_{n} (r_{n}))\, dr_{1} \dots dr_{n-1}
\biggr\Vert_{k}dr
 \\
&
\leq\int_{\R} \dots\int_{\R} \w(r)^{m}
 \pa \alpha_{s_{1}+\eta_{1}}(f_{1}(r_{1}))
\dots
\alpha_{s_{n-1}+\eta_{n-1}}(f_{n-1} (r_{n-1}))
\alpha_{s_{n}+\eta_{n}}(f^{(l)}_{n} (r_{n}))\pa_{k} dr_{1} \dots dr_{n},
 \endaligned
\tag ** $$
We bound the normed expression in the integrand. We use (*) and
 $\eta_{1}=0$ and  $\eta_{2}=r_{1}$
in the first step, and proceed analogously through the remaining steps.
$$\aligned
&\pa \alpha_{s_{1}+\eta_{1}}(  f_{1}(r_{1}))
\dots
\alpha_{s_{n-1} + \eta_{n-1}}(f_{n-1}(r_{n-1}))
\alpha_{s_{n}+\eta_{n}}(f^{(l)}_{n}(r_{n}))\pa_{k} \\
&
\leq
C\w(s_{1})^{d}
\pa f_{1}(r_{1}) \pa_{k} \,\pa \alpha_{s_{2}-s_{1}+r_{1}}
(f_{2}(r_{2}))\dots
\alpha_{s_{n-1} - s_{1} + \eta_{n-1}}(f_{n-1}(r_{n-1}))
\alpha_{s_{n}-s_{1} +\eta_{n}}
(f^{(l)}_{n}(r_{n})) \pa_{k} \, \\
&\leq
C^{2}\bigl(\w(s_{1})\,\w(s_{2}-s_{1})\, \w(r_{1}) \bigr)^{d}
\pa f_{1}(r_{1}) \pa_{k}\,\pa f_{2}(r_{2})\pa_{k}
\,\pa \alpha_{s_{3}-s_{2}+r_{2}}(f_{3}(r_{3}))\\
&\qquad  \qquad \qquad \qquad \qquad \qquad \dots\dots\dots
\alpha_{s_{n-1} - s_{2} + \eta_{n-1}}(f_{n-1}(r_{n-1}))
\alpha_{s_{n} - s_{2} + \eta_{n}}(f_{n}^{(l)}(r_{n}))
\pa_{k} \\
& \dots \qquad \dots \qquad \dots \\
&
\leq C^{n}\bigl(\w(s_{1})\dots \w(s_{n}-s_{n-1}) \bigr)^{d}
\bigl(\w(r_{1}) \dots \w(r_{n}) \bigr)^{d}
\,\pa f_{1}(r_{1}) \pa_{k}
\dots
\pa f_{n-1}(r_{n-1}) \pa_{k}\,
\pa
f^{(l)}_{n}(r_{n}) \pa_{k}. \endaligned
$$
Hence (**) is bounded by
$$C^{n}\bigl(\w(s_{1})\w(s_{2}-s_{1})\dots \w(s_{n}-s_{n-1}) \bigr)^{d}
\pa f_{1} \pa_{m+d, 0, k} \dots
\pa f_{n-1} \pa_{m+d, 0,k}\,
\pa f_{n} \pa_{m+d, l, k}.$$
Thus the increasing seminorms
$ \pa f \pa_{k}' = \max_{l\leq k} \int_{\R} \w(r)^{k}
 \pa f^{(l)}(r) \pa_{k} dr$ on $\sona$
satisfy the condition
$$ \pa \alpha_{s_{1}}(f_{1}) \dots \alpha_{s_{n}}(f_{n}) \pa_{k}'
\leq C^{n} \bigl(\w(s_{1})\w(s_{2}-s_{1})\dots \w(s_{n}-s_{n-1})
 \bigr)^{d} \pa f_{1} \pa_{k+d}' \dots \pa f_{n} \pa_{k+d}' $$
of Theorem 3.1.18 of \cite{\Scha}.  It follows that $\alpha$ is an
$m$-tempered action on $\sona$.
\par
Now let $\pa \quad \pa_{n}$ be
seminorms topologizing $\sona$ as in Definition 1.1.3.
Let $B_{n}$ be the completion of $\sona$
in $\pa \quad \pa_{n}$, and let $\alpha^{(n)}$
denote the extension of $\alpha$ to $B_{n}$.
Then $B_{n}$ is a Banach algebra and $\alpha^{(n)}$ is
a strongly continuous action  of $\R$ by algebra automorphisms.
Let $B_{n, m}$ denote the Banach algebra of $m$-times differentiable
vectors for the action $\alpha^{(n)}$ on $B_{n}$.  Then $B_{n,n}$ is
dense and spectrally invariant in $B_{n}$. (For example, this
follows from the
estimates in \cite{\Schb, Theorem 2.2}.)
Since $\sona = \varprojlim B_{n}$, we have
$$ \Cal S(\R, A^{\infty}, \alpha) =\sona^{\infty} =
\bigl( \varprojlim B_{n} \bigr)^{\infty} =
\varprojlim B_{n}^{\infty} =
\varprojlim_{n} \varprojlim_{m} B_{n, m} = \varprojlim B_{n, n}.  $$
A direct application of Lemma 1.1.9 (2) now tells us that the
inclusion $\Cal S(\R, A^{\infty}, \alpha) \hookrightarrow \sona$ is an
isomorphism on $K$-theory.
\qed
\enddemo
\proclaim{Corollary 1.3.2}  Let $A$ be a Fr\'echet algebra with
$m$-tempered, strongly continuous (not necessarily
smooth) action $\alpha$ of $\R$.  Then there is an isomorphism
$$ RK_{*}(A) \cong RK_{*+1}(\sona).$$
\endproclaim
\demo{Proof} The argument in the last paragraph of the proof of Theorem
1.3.1 shows that the inclusion $A^{\infty} \hookrightarrow A$ is an
isomorphism on $K$-theory.  We have the diagram
$$ \CD
RK_{*}(A^{\infty})    @>i_{*}>>   RK_{*}(A)\\
@V\theta VV         \\
RK_{*+1}(\Cal S(\R, A^{\infty}, \alpha)) @>i_{*}>> RK_{*+1}(\sona),
\endCD $$
where the bottom arrow is an isomorphism by Theorem 1.3.1, and $\theta$
is the Thom map from $\S 1$. The composition  of the three isomorphisms
gives the desired isomorphism on the right hand side.
\qed
\enddemo
\vskip\baselineskip
\subheading{ Definition 1.3.3}
We say that $\alpha$ acts {\it isometrically} on a Fr\'echet algebra
 $A$ if there exists a family $\bigl\{\pa \quad \pa_{k}\bigr\}$
 of submultiplicative seminorms topologizing $A$ such that
$\alpha$ acts isometrically for each $\pa \quad \pa_{k}$.
(This is equivalent to the existence of not necessarily
submultiplicative seminorms $\pa\quad\pa_k$
such that $\pa \alpha_{s}(a) \pa_{k}$
is a bounded function of $s\in \R$ for each $a \in A$ and $k$.
See the proof of Lemma 1.2.3 and Remark 1.2.4.)
Define
$$L_{1}(\R, A, \alpha) = \biggl\{\, f \colon \R \rightarrow A\,
 \biggl| \,f \text{ measurable and }
\int_{\R} \pa f(r) \pa_{k} dr < \infty\, \,\text{ for all }
k \in \N\, \biggr\}.$$
(See \cite{\Scha, \S 2} for details.)
\par
If $\alpha$ is an isometric action on $A$, then $L_{1}(\R, A, \alpha)$
is a Fr\'echet algebra under convolution. However, if $\alpha$ is not
isometric,  then $L_{1}(\R, A, \alpha)$ is unlikely to be an algebra.
Note that if $\alpha$ is isometric and $A$ is a Banach algebra,
then so is $L_{1}(\R, A, \alpha)$.
\proclaim{Theorem 1.3.4}  Let $A$ be a Fr\'echet algebra with
isometric action $\alpha$ of $\R$.  Then the inclusion map
$$ \sona \hookrightarrow L_{1}(\R, A, \alpha)$$
is an isomorphism on $K$-theory.
\endproclaim
\demo{Proof}
Let $\pa \quad \pa_{k}$ be submultiplicative $\alpha$-isometric
seminorms for $A$.  Let $A_{k}$ be the completion of
$A/{\text{Ker}(\pa \quad \pa_{k})}$
in $\pa \quad \pa_{k}$, and let $B_{k} = L_{1}(\R, A_{k}, \alpha)$.
Then the inverse limit $\varprojlim B_{k}$ is equal to
$L_{1}(\R, A, \alpha)$.  Define
$$ \pa f \pa'_{n,k} =
\sum_{i+j = n} \int_{\R} \w(r)^{i} \pa f^{(j)}(r) \pa_{k} dr,
\quad f \in \sona, $$
where $\w(r)= 1+|r|$.
Let $B_{n, k}$ be the completion of
$\sona/{\text{Ker}(\pa \quad \pa'_{n, k})}$
in $\pa \quad \pa'_{n, k}$.  The seminorms $\pa \quad \pa'_{n, k}$
are submultiplicative under convolution since the action is isometric.
Moreover, omitting the subscript $k$ everywhere for convenience, we have
$$
\aligned  \pa & f*g \pa'_{n}  = \sum_{i+j= n}
\int_{\R} \w(r)^{i} \,\biggl\Vert \, \int_{\R}
\,f(t) \,\alpha_{t}(g^{(j)}(r-t))\,dt \biggr\Vert dr\\
&
\leq \sum_{i+j= n} \int_{\R} \,\w(r)^{i} \,
 \int_{\R} \,\pa f(t)\pa \,\pa g^{(j)}(r-t))\pa \,dt  dr\\
&
\leq 2^{n} \sum_{i+j= n}\biggl(
\int_{\R} \int_{\R}
\, \w(t)^{i}\, \pa f(t) \pa \, \pa g^{(j)}(r-t) \pa\, dt dr
\\
&\qquad\qquad +
\int_{\R} \int_{\R}\,
  \pa f(t) \pa \,\w(r-t)^{i}\, \pa g^{(j)}(r-t) \pa\, dt dr
\biggr)\quad
\\
&\leq 2^{n} \biggl( \sum_{i_{1} + i_{2} + j_{1} + j_{2}=n}
\int_{\R}\int_{\R} \,\w(t)^{i_{1}} \pa f^{(i_{2})}(t) \pa
\,\w(r-t)^{j_{1}} \,\pa g^{(j_{2})}(r-t) \pa\, drdt
\biggr)\\&
\qquad\qquad\qquad\qquad\qquad
\qquad\qquad\qquad\qquad\qquad
\qquad\qquad\qquad \,\,\,
 + 2^{n}\pa f \pa_{0}^{\prime} \, \pa g \pa_{n}^{\prime}\\
& \leq
2^{n+1}\sum_{i+j= n}\, \pa f\pa'_{i} \, \pa g \pa'_{j}.\endaligned $$
This condition implies the  spectral invariance of $B_{n, k}$ in
$B_{k}$ for every $n$ \cite{\Schb, (1.11), Theorem 1.17}. (The main
point here is that the sum on the right hand side is over $i+j= n$.)
In particular, $B_{k,k}$ is spectral invariant in $B_{k}$.  We have
$$ \sona = \varprojlim_{k} \varprojlim_{n} B_{n, k}
= \varprojlim B_{k, k}, $$
so by Lemma 1.1.9(2), we have that
$\sona \hookrightarrow  L_{1}(\R, A, \alpha)$
is an isomorphism on $K$-theory.
\qed
\enddemo
\proclaim{Corollary 1.3.5}
Let $A$ be a  Fr\'echet algebra with isometric and strongly
continuous action $\alpha$ of $\R$.  Then there is an isomorphism
$$RK_{*}(A) \cong RK_{*+1}(L_{1}(\R, A, \alpha)).$$
If $A$ is a Banach algebra, then so is $L_{1}(\R, A, \alpha)$,
and  $RK$ may be replaced with $K$ in the isomorphism.
\endproclaim
\demo{Proof}
This is similar to Corollary 1.3.2.  By Theorem 1.3.4
and Corollary 1.3.2  we have isomorphisms
$$
\CD RK_{*}(A) @>\theta >>
RK_{*+1}(\sona) @>i_{*}>> RK_{*+1}(L_{1}(\R, A, \alpha)).
\endCD $$
Thus by taking the composition, we get the desired isomorphism.
For replacing $RK$ with $K$, apply \cite{\Ph, Corollary 7.8}.
\qed
\enddemo
\vskip\baselineskip
\heading \S 2 Smooth Crossed Products by $\Z$ \endheading
\par
We introduce the smooth mapping torus, and use it and the
smooth Thom isomorphism of \S 1 to
obtain a smooth version of the Pimsner-Voiculescu exact sequence.
This is an imitation of the proof of the Pimsner-Voiculescu exact
sequence for C*-algebras in \cite{\Bl, \S 10.3-4}.
\subheading{2.1 Definition}
Let $\alpha$ be an automorphism of a Fr\'echet algebra $A$.
Define the {\it mapping torus $M(A, \alpha)$} to be
$$ M(A, \alpha) = \biggl\{ \varf \colon [0, 1] \rightarrow A,
\, \varf \,\text{continuous } \biggl| \,\varf (1) = \alpha (
\varf (0)) \biggr\}.  $$
Then $M(A, \alpha)$ is a Fr\'echet algebra under pointwise
multiplication, with the topology of uniform convergence.
Define the {\it smooth mapping torus $M_{\infty}(A, \alpha)$}
to be the set of functions
$$ M_{\infty}(A, \alpha) = \biggl\{ \varf \colon [0, 1] \rightarrow A,
\, \varf \,\text{differentiable } \biggl| \,\varf^{(l)}(1) = \alpha (
\varf^{(l)}(0)),
\quad l \in \N\, \biggr\}.  $$
We topologize $M_{\infty}(A, \alpha)$ by uniform convergence of
derivatives.  Under pointwise multiplication, the sum of the sup norms
of the first $n$ derivatives is submultiplicative modulo constants.
Hence $M_{\infty}(A, \alpha)$ is a Fr\'echet algebra.
\par
We define the {\it smooth $\alpha$-suspension $S_{\infty}(A, \alpha)$ }
to be the closed subalgebra of $M_{\infty}(A, \alpha)$ given by
$$
S_{\infty}(A, \alpha) =
\biggl\{ \varf \colon [0, 1] \rightarrow A, \,
\varf \text{ differentiable }\biggl|
\, \varf(1)=\varf(0) = 0,\, \varf^{(l)}(1) =
\alpha(\varf^{(l)}(0)),\, l\in \N \, \biggr\}.  $$
\proclaim{2.2 Lemma} If $A$ is a  Fr\'echet algebra with automorphism
$\alpha$, and $SA = \{\, f \in C(\T, A)\, |\, \,f(0) = 0 \,\,\}$
is the continuous suspension of $A$, then the inclusion
$S_{\infty}(A, \alpha) \hookrightarrow SA$ is an isomorphism on
$K$-theory.  In particular, there is a natural isomorphism
$ RK_{*}(A) \cong RK_{*+1}(S_{\infty}(A, \alpha))$.
\endproclaim
\demo{Proof}
The proof is similar to that of Lemma 1.1.10(3). We use $[0,1]$
in place of $\R$, and replace the rapid decay conditions on the
derivatives with the condition that the derivatives match correctly at
$0$ and $1$, as in Definition 2.1. The Fourier transform argument
becomes unnecessary. Details are omitted.
\qed \enddemo
\proclaim{2.3 Lemma} Let $A$ be a Fr\'echet algebra, and let $\alpha$
be an automorphism of $A$. Then there is a natural six term exact
sequence:
$$ \CD RK_{1}(A) @>\delta_0 >> RK_{0}(M_{\infty}(A, \alpha))
@> >> RK_{0}(A) \\
@A {\text{id}} - (\alpha^{-1})_*
AA @. @VV {\text{id}} - (\alpha^{-1})_* V \\
RK_{1}(A) @< <<
RK_{1}(M_{\infty}(A, \alpha)) @< \delta_1 << RK_{0}(A). \\
\endCD $$
Here, $\delta_i$ is  the isomorphism of Lemma 2.2,
followed by the inclusion of $S_{\infty}(A, \alpha)$ in
$M_{\infty}(A, \alpha)$.
\endproclaim
\demo{Proof}
This is essentially Proposition 10.4.1 of \cite{\Bl}. We give a
different
proof, since the one in \cite{\Bl} uses a formula for the connecting
homomorphisms which has not been proved in representable K-theory.
\par
We will show that the sequence of the lemma is the one associated to the
short exact sequence
$$ \CD 0 @> >> S_{\infty}(A, \alpha) @> >>
 M_{\infty}(A, \alpha) @> {\text{ev}}_{1}>> A @> >> 0. \endCD \tag * $$
The only thing requiring proof is the identification of the vertical
maps ${\text{id}} -
(\alpha^{-1})_*$, and it suffices to do the one on the left.
\par
We replace $(*)$ by
$$ \CD 0 @> >> SA @> >>
 M(A, \alpha) @> {\text{ev}}_{1}>> A @> >> 0. \endCD \tag ** $$
The inclusion maps $S_{\infty}(A, \alpha) \rightarrow SA$ and
$A \rightarrow A$ are isomorphisms on $K$-theory. Therefore the Five
Lemma implies that $M_{\infty}(A, \alpha) \rightarrow M(A, \alpha)$
is too. This justifies the replacement.
\par
Let $\beta \colon RK_1 (A) \rightarrow RK_0 (SA)$ be Bott periodicity,
and let $\partial \colon RK_1 (A) \rightarrow RK_0 (SA)$
be the connecting homomorphism from $(**)$. We show that
$\beta^{-1} \circ \partial  = {\text{id}} - (\alpha^{-1})_*$.
Using the construction of the boundary map in  Proposition 3.11
and Remark 3.12 of \cite{\Ph},
we identify $\beta^{-1} \circ \partial $ with the composition
$$ \CD RK_1 (A) @> \beta >> RK_0 (SA) @> \varphi_{*}^{-1}\circ \nu_*  >>
         RK_0 (SA) @> \beta^{-1} >> RK_1 (A). \endCD $$
Following Remark 3.12 of \cite{\Ph}, and being careful of the
conflicting notation, we obtain $\varphi \colon SA \rightarrow D$ and
$\nu \colon SA \rightarrow D$ as:
$$ D = \biggr\{ (g,f)\in C([0,1], A) \oplus C([0,1], A) \biggr| \quad
   g(1)= \alpha (g(0)), f(0)=g(1), \text{and}\  f(1) = 0 \biggl\}, $$
$$ \varphi (g) = (g,0), \quad \quad \quad \quad \text{and} \quad
 \quad \quad \quad \nu (f) = (0, f). $$
Then $\nu$ is homotopic to the map
$\mu \colon SA \rightarrow D$ given by $\mu(f) = (g,0)$, with
$g(t/2) = \alpha^{-1}(f(1-t))$ and
$g((t+1)/2) = f(t)$ for $t \in [0,1]$.
A homotopy is given by  $\nu_{s} (f) = (g_{s}, h_{s})$, with:
$$ h_{s}(t) = \cases
f(t + s) & \qquad t+ s \leq 1
\\ 0 & \qquad t + s >1 \endcases $$
$$ g_{s}(t/2) = \cases \alpha^{-1}(f(s-t))  & \qquad t + (1-s) \leq 1
 \\
 0  & \qquad t+ (1-s) > 1 \endcases $$
$$ g_{s}(t/2 + 1/2) = \cases f(s +t-1)  & \qquad (1-t) + (1-s) \leq 1 \\
 0  & \qquad (1-t)+ (1-s) > 1 \endcases $$
where $s \in [0,1]$ and $t \in [0,1]$.
So
$\varphi_{*}^{-1}\circ \nu_*  =  \varphi_{*}^{-1}\circ \mu_* $. Since
$\varphi$ is an isomorphism onto its image, and since $\varphi_*$ is an
isomorphism, we can replace $D$ by $SA$ and
$ \varphi_{*}^{-1}\circ \mu_* $
by $( \varphi^{-1}\circ \mu )_*$. This map is
$ - (S\alpha^{-1})_* + {\text{id}}$
by Lemma 6.2 of \cite{\Ph}. Therefore $\beta^{-1} \circ \partial
=  {\text{id}}- (\alpha^{-1})_*$
by naturality of $\beta$.
\qed \enddemo
\par
If $\alpha$ is a continuous action of the circle group $\T$ on  $A$
(automatically tempered and $m$-tempered since $\T$ is compact - see
Definition 1.1.3 and \cite{\Scha, Definition 2.2.4}), we define the
Fr\'echet algebra $\Cal S(\T, A, \alpha)$ to be
$C^{\infty}(\T)\otimes A$ with the twisted convolution product.  The
following proposition and proof also
hold for non $m$-convex Fr\'echet algebras $A$ and tempered, but
not necessarily $m$-tempered, actions $\alpha$ of $\R$ on $A$.
\proclaim {2.4 Proposition}
Let $\R$ act on a Fr\'echet algebra $A$ with $m$-tempered action
$\gamma$, and assume that the restriction of $\gamma $ to $\Z$ is the
trivial action.  Then $\gamma$ factors to an action $\overline \gamma$
of $\T= \R/\Z$ on $A$, and $\beta= (\widehat {\overline \gamma})_1$
is an automorphism of $\Cal S(\T , A, {\overline \gamma})$.
The smooth crossed product $\Cal S(\R, A, \gamma)$ is isomorphic to the
Fr\'echet algebra $M_{\infty}
(\Cal S(\T,A, {\overline \gamma}), {\beta})$.
\endproclaim
\demo{Proof}
The first statement is trivial. In the proof of the second,
we use the convention $e(t) = e^{2\pi i t}$.
Let $B=\Cal S(\R,  A, \gamma)$ and $C=M_{\infty}(
\Cal S(\T, A, {\overline \gamma}), { \beta})$.
We define a map $\Psi \colon B \rightarrow C$ by
$$ \Psi(\varf)(s) (t + \Z) = \sum_{k \in \Z}
\,e(s(t+k))\,\, \varf(t + k), \quad s\in [0, 1],\quad t \in \R.  $$
Since $\varf$ is a Schwartz function, it follows that
$\Psi(\varf)(s)(t+\Z)$ is infinitely differentiable in both variables.
Also, the functions $t \mapsto \Psi(\varf)^{(l)}(s)(t+\Z)$ and their
derivatives are periodic in $t$, so $\Psi(\varf)^{(l)}(s)$ does
in fact define a $C^{\infty}$ function from $\T$ to $A$ for each
$s \in [0, 1]$ and $l \in \N$. (We denote the $l$-th derivative in
$s$ by $\Psi(\varf)^{(l)}(s)(t+\Z)$ and the $l$-th derivative in $t$ by
$\Psi(\varf)(s)^{(l)}(t + \Z)$.)
A quick calculation shows that $\Psi(\varf)^{(l)}(1)(t + \Z) =
e(t) \Psi(\varf)^{(l)}(0)(t+\Z)={ \beta}(\Psi(\varf)^{(l)}(0))(t+\Z)$.
It follows that $\Psi(\varf) \in C$.  Moreover, it is not
to hard to check that if $\varf_{n} $ converges to $\varf$
in the Schwartz topology of $B$, then the
derivatives of $\Psi(\varf_{n})$ converge uniformly to the
derivatives of $\Psi(\varf)$.  So $\Psi$ is a continuous
linear map of Fr\'echet spaces.
\par
We verify that $\Psi$ is an algebra homomorphism.
$$ \aligned
\Psi(\varf * g)(s)(t + \Z)
%
%
& =  \sum_{k \in \Z}\,e(s(t+k))\,\int_{\R} \varf(r) \,
\gamma_{r}(g(t+k-r)) \,dr \\
%
%
& =  \sum_{k, l \in \Z} \int_{0}^{1}\, e(s(w+l))\,\, \varf(w+ l)\,\,\,
e(s(t+k-w))\,\,{\gamma}_{w}(g(t+k -w))\,dw\\
&
 = \biggl(\Psi(\varf)*\Psi(g)\biggr)(s)(t + \Z).
\endaligned
$$
\par
We construct an inverse for $\Psi$.  Define a
map $\Phi \colon C \rightarrow B$ by
$$ \Phi(\xi)(r) = \int_{0}^{1} \,e(-wr)\,\xi(w)(r+\Z)\, dw, \quad
r \in \R.
$$
We show that $\Phi (\xi)$ lies in $ B$.
First note that by integration by parts, we have
$$
(2\pi i r)^{p}
\int_{0}^{1} \,e(-wr)\,\xi(w)(r+\Z)\, dw
=
\int_{0}^{1} \,e(-wr)\,\xi^{(p)}(w)(r+\Z)\, dw. \tag *$$
Here the boundary condition
$\xi^{(k)}(1) (r+\Z) = e(r)\xi^{(k)}(0) (r+\Z)$ is essential for
eliminating terms at the endpoints. From (*), it follows that
$$  \biggl\Vert
\biggl({d\over{dr}}\biggr)^{l}
\,r^{p}\,\Phi(\xi)(r)\biggr\Vert_{m} \leq
\biggl({1\over{2\pi}}\biggr)^{p}\int_{0}^{1} \, \biggl\Vert
\biggl({d\over{dr}}\biggr)^{l}
e(-wr)\,\xi^{(p)}(w) (r+\Z)\biggr\Vert_{m} \, dw
 $$
for each continuous seminorm $\pa \quad \pa_{m}$ on $A$.
Since the integrand is a bounded continuous function, we see that
$\sup_{r\in\R}\pa \biggl({d\over{dr}}\biggr)^{l} \,r^{p}\,
\Phi(\xi) (r)\pa_{ m} <\infty$
for any $l, p,m \in \N$.  Using the commutation relation
$[{d\over{dr}} ,r] = 1$ repeatedly, we see that
the seminorms $\pa \Phi(\xi)\pa_{p,l,m}
= \sup_{r\in \R}\pa r^{p} \, \Phi(\xi)^{(l)}(r)\pa_{m}$
on the smooth crossed product $B$ are all
finite.  Hence $\Phi(\xi) \in  B$.
\par
For $\xi \in C$, we have
$$ \aligned \Psi \circ \Phi(\xi) (s)(t + \Z) &=
 \sum_{k \in \Z}\int_{0}^{1} e(s(t+k)) \,e(-w(t + k))\,
\xi(w)(t+k+\Z)\, dw \\
& =  e(st) \, \sum_{k \in \Z} e(sk) \, \int_{0}^{1} e(-wk) \,
\bigl[e(-wt)\xi(w)(t+\Z)\bigr]\, dw \\
& = e(st) \bigl[e(-st)\xi(s)(t + \Z)\bigr] = \xi(s)(t + \Z).
\endaligned $$
The third step uses the fact that the term in brackets is a smooth
periodic function of $w$, and is therefore the sum of its Fourier
series.
Similarly, one checks that  $\Phi \circ \Psi(\varf) = \varf$ for
$\varf \in  B$. Thus, by the open mapping theorem, $\Psi$ is an
isomorphism of Fr\'echet algebras.
\qed \enddemo
\par
For an $m$-tempered action $\alpha$ of $\Z$ on a Fr\'echet algebra $A$,
we define $\Cal S(\Z, A, \alpha)$ to be $\Cal S(\Z) \otimes A$
with twisted convolution  multiplication.  Then
$\Cal S(\Z, A, \alpha)$ is a Fr\'echet algebra, by
\cite{\Scha, Theorem 3.1.7}.  We
define the dual action $\widehat \alpha$ of $\T$ on
$\Cal S(\Z, A, \alpha)$ by ${\widehat \alpha}_{z}(f)(n) = e^{2\pi i zn}
f(n)$.
\proclaim{2.5 Lemma }
There is a natural isomorphism (Takesaki-Takai duality)
$$ d_{\alpha}\colon
\Cal S(\T, \Cal S(\Z, A, \alpha), {\widehat \alpha})
\rightarrow \Kinf \otimes A. $$
This isomorphism is equivariant for the second dual action
${\widehat {\widehat {\alpha}}}_{n} (f)(z,m) =
e^{2 \pi i nz} f(z,m)$ on the left
and the action on the right generated by the automorphism
$e_{jk} \otimes a  \mapsto e_{j-1, k-1} \otimes \alpha_1 (a)$.
\endproclaim
\demo{Proof}
The proof is essentially the same as for $\R$ in place of $\Z$.
We identify $\Cal S(\T, \Cal S(\Z, A, \alpha), {\widehat \alpha})$
with the Fr\'echet space $C^{\infty}(\T) \otimes \Cal S(\Z) \otimes A$.
Then we define $d_{\alpha}$ to be the composite $\Gamma \circ \Cal F$,
where $\Cal F$ is the Fourier transform in the first variable,
$$ \Cal F(f)(k, n) = \int_{\T} e^{2\pi i uk} \,f(u,n)\, du, $$
for $n,k \in \Z$, and $\Gamma (\eta)(k, n) = \alpha_{-k}(\eta(n, k-n))$.
One checks that $d_{\alpha}$ transforms the multiplication on
$C^{\infty}(\T) \otimes \Cal S(\Z) \otimes A$ into the usual matrix
multiplication on $\Cal S(\Z) \otimes \Cal S(\Z) \otimes A \cong
\Kinf \otimes A$, and transforms the second dual action as described
in the lemma.
\qed \enddemo
\proclaim{2.6 Theorem (Pimsner-Voiculescu Exact Sequence)}
Let $\alpha$ be an $m$-tempered action of $\Z$ on a Fr\'echet algebra
$A$.  Then we have the following natural six term exact sequence in
$K$-theory.
$$ \CD RK_{0}(A) @>{\text{id}}-(\alpha_{-1})_{*} >> RK_{0}(A)
@> \iota_{*}>> RK_{0}(\Cal S(\Z, A, \alpha)) \\
@A\partial AA @. @VV \partial V \\
RK_{1}(\Cal S(\Z, A, \alpha)) @< \iota_{*}<<
RK_{1}(A) @< {\text{id}}-(\alpha_{-1})_{*} << RK_{1}(A), \\
\endCD $$
where $\iota \colon A \rightarrow \Cal S(\Z, A, \alpha)$ is the
inclusion.
\endproclaim
\demo{Proof}
Let $ {\overline \gamma} = {\widehat \alpha}$ be the dual action of
$\T$ on $B= \Cal S(\Z, A, \alpha)$.  We lift ${\overline \gamma}$
to an action $\gamma$ of $\R$ on $B$. Apply Proposition 2.4 to this
action, and use Lemmas 2.3 and 2.5 to compute the $K$-theory of the
resulting mapping cylinder. We get an exact sequence

$$ \CD RK_{1}(\Kinf \otimes A) @> >> RK_{0}(\Cal S(\R, B, \gamma))
@> >> RK_{0}(\Kinf \otimes A) \\
@A {\text{id}} - (\lambda \otimes \alpha_{-1})_*
AA @. @VV {\text{id}} -
(\lambda \otimes \alpha_{-1})_*  V \\
RK_{1}(\Kinf \otimes  A) @< <<
RK_{1}(\Cal S(\R, B, \gamma)) @<  << RK_{0}(\Kinf \otimes A), \\
\endCD $$
where $\lambda \in {\text {Aut}}(\Kinf)$ is conjugation by the
bilateral shift.
Since $\lambda$ is homotopic to the trivial automorphism
of $\Kinf$, we can use
stability in $K$-theory to delete it and $\Kinf$. The smooth Thom
isomorphism allows us to replace $RK_{i}(\Cal S(\R, B, \gamma))$ by
$RK_{1-i}(B)$. Rotating the exact sequence two spaces counterclockwise
then gives the exact sequence of the theorem, except for the
identification of the maps $\iota_*$.
\par
The proof of this identification is not given in \cite{\Bl}, so we give
details here. We must show commutativity of the following diagram, with
maps as described below:
$$ \CD RK_0(A) @> >> RK_0 (\Kinf \otimes A) @> >>
                RK_1 ( S_{\infty}(\Kinf \otimes A, \beta)) \\
@V i_{*} VV @. @VV i_{*} V \\
RK_0 (B) @> >> RK_1 (\Cal S(\R, B, \gamma)) @>  >>
                 RK_1 ( M_{\infty}(\Kinf \otimes A, \beta)). \\
\endCD $$
The vertical maps are induced by the inclusions. The maps across the top
are stability followed by Bott periodicity (Lemma 2.2).
The maps across the bottom
are the Thom isomorphism (Theorem 1.2.7) and the map induced by the
isomorphisms in Proposition 2.4 and Lemma 2.5.
\par
We consider an arbitrary class in $RK_0 (A)$. Replacing $A$ by
$M_2 ((\Kinf \otimes A)^+)$ throughout, we may assume our class is
represented by an idempotent $p \in A$. (Compare with the argument using
the second diagram in the proof of Theorem 1.2.7.) We begin by computing
the image of $[p]$ via the lower left corner.
\par
Let $q$ be the image of $p$ in $B$.
For the next step, choose  $f_0 \in \Cal S(\R)$ whose Fourier
transform $\widehat \varf_0$ has support in $(0,1)$ and such that
$1 + \widehat \varf_0$ is an invertible element of $\Cal S(\R)^+$ (for
pointwise multiplication in $\Cal S(\R)$) which represents the standard
generator of $RK_1 (\Cal S(\R))$. Regarding elements of
$\Cal S(\R, B, \gamma)$ as functions from $\R \times \Z$ to $A$, let
$\varf \in \Cal S(\R, B, \gamma)$ be
$\varf (s,n) = \delta_{0n} \varf_0 (s)p$, that is,
$\varf = \varf_0 \otimes q$. Since $q$ is $\gamma$-invariant, this
formula actually yields an invertible element $1+\varf$ in the
unitization of each crossed product $ \Cal S(\R, B, \gamma^{(r)})$,
where  $\gamma_{s}^{(r)} = \gamma_{rs}$ for  $r \in [0,1]$.
Definition 1.2.1 (of the Thom map)
shows that $1+\varf$ represents the image of $[q]$ in
$RK_1 (\Cal S(\R, B, \gamma))$.
\par
We now compute the image of $\varf$ under the isomorphisms of
Proposition 2.4 and Lemma 2.5. We first apply the map $\Psi$ from the
proof of Proposition 2.4, and then apply pointwise on $[0,1]$ the maps
$\Cal F$ and $\Gamma$ from the proof of Lemma 2.5. Writing the result as
a function of $s \in [0,1]$ and $l,n \in \Z$, combining exponentials,
and rearranging slightly, we get
$$ (\Gamma \circ \Cal F \circ \Psi)(\varf)(s,l,n) =
   \sum_{k \in \Z} \int_{0}^{1} {\text{exp}}(2 \pi i(tn+st+sk)) \,
    \alpha_{-l}(\varf(t+k, l-n)) \, dt. $$
Since $nk$ is an integer, we can rewrite the exponential as
${\text{exp}}
(2\pi i (t+k)(s+n))$. Putting in the definition of $\varf$, we
find that $(\Gamma \circ \Cal F \circ \Psi)(\varf)(s,l,n) = 0$ for
$l \neq n$, and
$$ \aligned
 (\Gamma \circ \Cal F \circ \Psi)(\varf)(s,n,n) &=
  \alpha_{-n}(p) \sum_{k \in \Z} \int_{0}^{1} {\text{exp}}
(2 \pi i (t+k)(s+n))
                                \,  \varf_0(t+k) \, dt \\
   &=\alpha_{-n}(p) \widehat \varf_0 (s+n).
\endaligned $$
Since $s \in [0,1]$, we have $\widehat \varf_0 (s+n) = 0$  for
$n \neq 0$. We conclude that the image of $[p]$ in
$RK_1 ( M_{\infty}(\Kinf \otimes A, \beta))$ is represented by
$a = 1 + \bigl( \widehat \varf_0 \bigr|_{[0,1]} \bigl) \otimes e_{00}
   \otimes p$,
where $e_{00} \in \Kinf$ is $e_{00}(l,n) = \delta_{0l}\delta_{0n}$.
\par
We now go via the upper right corner. The image of $[p]$ in
$RK_0 (\Kinf \otimes A)$ is $[e_{00} \otimes p]$. One easily checks that
the image of this under Bott periodicity can be obtained as
$$ \biggl[1 + \bigl( \widehat \varf_0 \bigr|_{[0,1]} \bigl) \otimes
           e_{00} \otimes p \biggl]  \in
      RK_1 ( S_{\infty}(\Kinf \otimes A, \beta)), $$
by the choice of $\widehat \varf_0$. But this is just $a$. So we have
shown that the diagram commutes.
\qed
\enddemo
\vskip\baselineskip
We define the notion of an isometric action of $\Z$ on
a Fr\'echet algebra just as we did for an action of $\R$ in
Definition 1.3.3.
Let
$$L_{1}(\Z, A, \alpha) = \biggl\{\, f \colon \Z \rightarrow A\,
 \biggl\vert \, \sum_{m \in\Z} \pa f(m) \pa_{k}  < \infty\,\,
 {\text{ for all }} k \in \N\, \biggr\}.$$
If $\alpha$ acts isometrically on $A$, then $L_{1}(\Z, A, \alpha)$ is
a Fr\'echet algebra under convolution.
\proclaim{2.7 Lemma}  Let $A$ be a Fr\'echet algebra with
isometric action $\alpha$ of $\Z$.
Then the inclusion map
$$ \Cal S(\Z, A, \alpha) \hookrightarrow L_{1}(\Z, A, \alpha)$$
is an isomorphism on $K$-theory.
\endproclaim
\demo{Proof}
Let $\pa \quad \pa_{k}$ be increasing
submultiplicative $\alpha$-isometric
seminorms for $A$.  Let $A_{k}$ be the completion of
$A/{\text{Ker}(\pa \quad \pa_{k})}$
in $\pa \quad \pa_{k}$, and let $B_{k} = L_{1}(\Z, A_{k}, \alpha)$.
Then the inverse limit $\varprojlim B_{k}$ is equal to
$L_{1}(\Z, A, \alpha)$.  Define
$$ \pa f \pa_{n,k} =
 \sum_{m \in \Z} (1+|m|)^{n}\, \pa f(m) \pa_{k},
\quad f \in \Cal S(\Z, A, \alpha). $$
Let $B_{n, k}$ be the completion of
$\Cal S(\Z, A, \alpha)/{\text{Ker}(\pa \quad \pa_{n, k})}$
in $\pa \quad \pa_{n, k}$.
The seminorms $\pa \quad \pa_{n, k}$
are submultiplicative under convolution since the action is isometric.
Moreover, we have
$\pa f*g \pa_{n, k} \leq 2^n \bigl(\pa f \pa_{n,k} \pa g \pa_{0,k}
+ \pa f \pa_{0, k} \pa g \pa_{n, k}\bigr)$.
The proof is similar to, but simpler than, the computation in the proof
of Theorem 1.3.4. Now finish just as there.
\qed
\enddemo
\proclaim{2.8 Corollary} If $\alpha$ is an isometric action of $\Z$
on a  Fr\'echet algebra $A$, then we have the following natural
six term exact sequence in $K$-theory.
$$ \CD RK_{0}(A) @>{\text{id}}-(\alpha_{-1})_{*} >> RK_{0}(A)
@> \iota_{*}>> RK_{0}(L_{1}(\Z, A, \alpha)) \\
@A\partial AA @. @VV \partial V \\
RK_{1}( L_{1}(\Z, A, \alpha)) @< \iota_{*}<<
RK_{1}(A) @< {\text{id}}-(\alpha_{-1})_{*} << RK_{1}(A). \\
\endCD $$
If $A$ is a Banach
algebra, then so is $L_{1}(\Z, A, \alpha)$, and we may replace
$RK$ by $K$ everywhere in the diagram.
\endproclaim
\demo{Proof} This follows from Theorem 2.6 and Lemma 2.7.
For the statement about Banach algebras, see \cite{\Ph, Corollary 7.8}.
\qed
\enddemo
\heading \S 3 Examples and Applications \endheading
\par
We give several examples and applications of our theorems.
\subheading{3.1 Example}
For  $G= \Z$ or $\R$, there are many examples of  actions $\alpha$
on locally compact spaces $M$ which induce an $m$-tempered
action on a suitable Fr\'echet algebra
${\Cal S} (M)$ of \lq\lq smooth\rq\rq
functions on $M$. (See \cite{\Schb, Examples 7.20,
6.26-7, 2.6-7} and \cite{\Scha, Example 5.18}.  The examples
\cite{\Scha, Example 5.19, 5.23-4, 5.26-7} also include
examples of $m$-tempered actions, along with many actions that
are {\it not} $m$-tempered  with respect to the absolute value
function - see Definition 1.1.3 and \cite{\Scha, Definition 3.1.1}.)
Our smooth Thom isomorphism and Pimsner-Voiculescu exact sequence apply
to these examples. We note, however, that the $K$-theory of the smooth
crossed product can often be computed directly from Lemma 1.1.9(1) and
the C*-algebra results, because $\Cal S(G, \Cal S(M), \alpha)$
is often   spectral invariant in  $C^{*}(G, C_{0}(M), \alpha)$.
(See \cite{\Schb, Corollary 7.16}.)
\par
Two specific examples are the canonical smooth subalgebra
$\Cal S(\Z, C^{\infty}(\T), \alpha)$
of the (rational or irrational) rotation algebra,
and $\Cal S(\R, \Cal S(H/K), \alpha)$, where $H$ is
the unipotent $3\times 3$ upper triangular real
matrices, $\R$ the subgroup corresponding to the
first row and second column, and $K$ any closed subgroup of $H$.
\subheading{3.2 Example}
Let  $G= \Z$ or $\R$,  and let $\varphi \colon A \rightarrow B$
be a homomorphism which is equivariant for the $m$-tempered actions
$\alpha$ and $\beta$ on the Fr\'echet algebras $A$ and $B$. Assume
that
$\varphi_* \colon RK_{*}(A) \rightarrow RK_{*}(B)$ is an isomorphism.
Then
$RK_{*}(\Cal S(G, A, \alpha)) \rightarrow RK_{*}(\Cal S(G, B, \beta))$
is also an isomorphism. If $\beta$ is isometric and $B$ is a Banach
algebra, then
$RK_{*}(\Cal S(G, A, \alpha)) \rightarrow K_{*}(L_1(G, B, \beta))$
is an isomorphism, and if $B$ is a C*-algebra, then
$RK_{*}(\Cal S(G, A, \alpha)) \rightarrow K_{*}(C^{*}(G, B, \beta))$
is an isomorphism. For $G = \R$, these statements are immediate from
the Thom isomorphisms for the various crossed products, and for
$G = \Z$ they follow from the Pimsner-Voiculescu exact sequences
via the Five Lemma.
\par
This argument applies without any assumption on $\varphi$ beyond that it
is a continuous equivariant homomorphism which is an isomorphism on
$K$-theory. But we are of course most intersted in the case where
$\varphi$ is the inclusion of a dense subalgebra $A$ of a Banach or
C*-algebra $B$, with $A$ Fr\'echet in its own topology. If $A$ is
spectral invariant in $B$, then $\varphi$ is automatically an
isomorphism on $K$-theory, by Lemma 1.1.9(1). It follows that the
inclusion $\Cal S(G, A, \alpha) \rightarrow L_1(G, B, \beta)$
(or $\Cal S(G, A, \alpha) \rightarrow C^{*}(G, B, \beta)$ )
is an isomorphism on $K$-theory. Cases in which the smooth crossed
product is spectral invariant in the C* crossed product are mentioned
in Example 3.1, but here the point is that one gets the the isomorphism
of the $K$-theory of the crossed products without knowing anything
about spectral invariance of the crossed products.
\par
For example, consider Theorem 12.5 of
\cite{\Ns}, and assume that the action is $m$-tempered, so that Nest's
crossed product is the same as ours. Using $RK_*$ in place of Nest's
$K_{*}^{p}$ for the dense subalgebras, we get that if the inclusion of
the dense subalgebra is an isomorphism on $K$-theory
before taking crossed products,
then it is again an isomorphism (not just surjective) after taking
crossed products.
\par
The discussion above generalizes to  a closed subgroup $G$ of
a connected, simply connected, nilpotent Lie group  $H$.  Such a group
$G$ is a semidirect product $Z\rtimes Z \rtimes \dots \rtimes
\Z \rtimes \R \rtimes \R \rtimes \dots\rtimes \R$.
Let $\alpha$ and $\beta$ be an $m$-tempered
actions of $G$ on $A$ and $B$ for which $\varphi$ is equivariant.
(Here $m$-temperedness is meant with respect to the gauge from
\cite{\Scha, Theorem 1.5.13}.  See  also \cite{\Scha, Example
1.5.14} and \cite{\Scha, Definition 3.1.1}.)
Then, decomposing the smooth
crossed products $\Cal S(G, A, \alpha) $ and
$\Cal S(G, B, \beta)$ (defined in \cite{\Scha, \S 2})
into successive smooth crossed products  by $\R$
and then by $\Z$, we see by iterating the above results that
$ \Cal S(G, A, \alpha) \rightarrow \Cal S(G, B, \beta)$
is an isomorphism on $K$-theory.
\par
Consider the special case that $B$ is a C*-algebra, and $\beta$
is an action via *-automorphisms.  Then the inclusion map
$\Cal S(G, B, \beta) \hookrightarrow C^{*}(G, B, \beta)$
is spectral invariant with dense image
by \cite{\Schb, Corollary 7.16}. (The same result
in the special case $G=\Z$ is also done in \cite{\Bo, Theorem 2.3.3}.)
Hence $\Cal S(G, A, \alpha) \rightarrow C^{*}(G, B, \beta)$
is an isomorphism on $K$-theory.
\subheading{3.3 Example}
Let $U$ be the annulus
$\{\,\,z \in \C \quad | \quad 1/2 < |z| < 2\,\,\}$.
Let $A(U)$ be the Banach algebra of continuous functions from
$\overline U$ to $\C$ which are holomorphic on $U$. Let $C(U)$ be
the Fr\'echet algebra of all continuous functions on $U$, with the
topology of uniform convergence on compact subsets, and let $H(U)$
be the closed subalgebra of holomorphic functions on $U$. Then the
rotation by an irrational angle $\theta$ defines an isometric action
on each of these algebras. The maximal ideal spaces of these algebras
are either $\overline U$ or $U$, and the inclusion of the unit circle
$\T$ into each of these spaces is a homotopy equivalence. Using Theorem
7.15 of \cite{\Ph}, it is easy to see that the restriction map from
any of these algebras to $C(\T)$ is an isomorphism on $K$-theory.
It follows from the previous example that $L_1 (\Z, A(U), \theta)$,
 $L_1 (\Z, C(U), \theta)$, and  $L_1 (\Z, H(U), \theta)$ all have
the same $K$-theory as the irrational rotation C*-algebra $A_{\theta}$.
The same holds for the smooth crossed products
$\Cal S (\Z, A(U), \theta)$, etc.
\par
Let A be the Banach algebra of continuous functions on $\T$ with
absolutely convergent Fourier series, with the norm
$\pa \varf \pa = \sum_{n \in \Z} | \widehat \varf (n)|$.
A similar argument to the above also shows that
 $L_1 (\Z, A, \theta) \rightarrow A_{\theta}$
is an isomorphism on $K$-theory.
\par
In place of $U$ above, one can also use $\C$, the open unit disk $D$,
or $\C - \overline D$. The $K$-theory of the crossed product is always
the same as that of $C^{*}(\Z)$ for $\C$ and for $D$, and the same as
that of $A_{\theta}$ for $\C - \overline D$. One can also let $\R$
act on all these algebras via rotation, with similar results, or one
can let $\R$ act by translation on similar algebras gotten by using
a half plane or strip in place of a disk or an annulus.
\par
In some of the cases above, the $K$-theory of the crossed product
can be computed by other methods. For example, the disk algebra $A(D)$
is equivariantly homotopy equivalent to $\C$, via evaluation at $0$.
Therefore,  $L_1 (\Z, A(D), \theta)$ is homotopy equivalent to
$L_1 (\Z)$. Also, $\Cal S (\Z, C^{\infty} (\T), \theta)$ is dense and
spectral invariant in both $L_1 (\Z, A, \theta)$ and $A_{\theta}$
\cite{\Schb, Theorem 6.7, Corollary 7.16},
so that  $L_1 (\Z, A, \theta) \rightarrow A_{\theta}$ is an isomorphism
on $K$-theory by Lemma 1.1.9(1). However, the methods vary from case to
case, and it is not obvious that they apply to all the cases above.
\par
For other examples, consider the irrational rotation on $H^{\infty}(D)$,
$H^{\infty}(U)$, or the measure algebra $M(\T)$. We will not attempt to
compute the $K$-theory of
these crossed products. (Note that the action of $\R$ on these algebras
is not continuous.)
\subheading{3.4 Example}
We give an example of a smooth crossed product of $\R$ with $A$
in which $A$ is a subalgebra of a Banach algebra but not
spectrally invariant in any  Banach algebra.
The algebra in this example is thus essentially different from the
algebras of unbounded functions in Example 3.3.
For $n \in \N$, define
$$ A_n = \biggl\{ \, f \in L_{1}(\R) \, \biggl| \,
\int_{\R} e^{n |r|} |f(r) | dr < \infty\, \biggr\}
\quad \quad \quad \quad \quad \text{and} \quad \quad \quad \quad \quad
A = \bigcap_{n \in \N} A_n = \varprojlim A_{n}.$$
Then each $A_n$ is a Banach algebra under convolution and
$A$ is a Fr\'echet algebra.
Let $\R$ act on $A$  via $\alpha_{s}(f)(r) = e^{i sr} f(r)$.
This action  is isometric and hence $m$-tempered. We will compute the
$K$-theory of the smooth crossed product, a somewhat strange subalgebra
of $\Kinf$, by computing the $K$-theory of $A$ and applying the Thom
isomorphism. (It is possible to compute the $K$-theory of this crossed
product more directly, but the computation is more complicated than for
$A$.)
\par
We compute $RK_{*} (A)$  using Theorem 7.15 of \cite{\Ph}, which
requires that we know the maximal ideal space $\text{Max}(A^{+})$
and its compactly generated topology.
Every complex number $z$ defines a character
$\eta_{z}(f) = \int_{\R} e^{ i z r} f(r) dr $ of $A$ (with obvious
extension to a character of the unitization $A^{+}$).
Note that  $A^{+}$ also has the
character $\eta_{\infty}(\lambda 1 + f )= \lambda$
which vanishes on $A$.
We show that these are in fact all the characters of $A^{+}$.
Let $\xi $ be any character of $A^{+}$.  If
$\xi $ vanishes on $A$, then clearly $\xi = \eta_{\infty}$.  Otherwise,
$\xi$ is a continuous linear map from
$A$ to the Banach space $\C$. Therefore there must be some $n$ such
that $\xi$ lifts to a continuous linear map from $A_{n}$ to $\C$.
It follows that there is $\chi_{0} \in L_{\infty}(\R)$ such that
$\xi(f) = \int_{\R} f(r) e^{n|r|} \chi_{0}(r) dr$.
Using $\xi(f*g) = \xi(f) \xi(g)$, a standard argument shows that
$e^{n|r|}\chi_{0}(r)$ coincides almost everywhere with
a continuous function ${\chi}$
which satisfies ${ \chi}(r+s) = {\chi}(r) {\chi}(s)$.
(See, for example, \cite{\Con, proof of Theorem VII.9.6}.)
It follows that ${ \chi}(r) = e^{ i zr} $ for some complex number
$z$.
Since $|\chi(r)| \leq \pa\chi_{0} \pa_{\infty} e^{n|r|}$,
we get $|\text{Im} (z)| \leq n$.  This shows that $\xi = \eta_{z}$
with $|\text{Im}(z)|\leq n$.
Thus, the maximal ideal space $\text{Max}(A^{+})$ is exactly
$\C \cup \{ \infty \}$.
\par
We further see from this argument that
$\text{Max}(A_{n}^{+})$ is  the set
$$ K_{n} = \{\,\, z \in \C\,\,|\,\, |\text{Im}(z)|\leq n \, \,\}\,\,
\cup \,\, \{\,\infty\,\}, $$
topologized as the one point compactification of the strip.
(We check that this is the right topology. By compactness and
metrizability, it suffices to show that if $z_k \rightarrow z$,
then $\eta_{z_k} \rightarrow \eta_z$. For $z \neq \infty$, this follows from
the Lebesgue Dominated Convergence Theorem.
If $z = \infty$, we use the following argument.
For $s \in [-n,n]$ define $f_{s}(r) = e^{-sr} f(r)$ for $f \in A_{n}$.
Then $s \mapsto {\widehat f_{s}}$ is a continuous
map from $[-n, n]$ to $C_{0}(\R)$ and so has compact
image.  Hence for each $\epsilon >0$
the image is contained in one of the
open sets $U_{\epsilon, m} = \{ \,\, \xi \in C_{0}(\R) \,\,| \,\,
|\xi(r) | < \epsilon {\text{ for all }} r \notin [-m, m] \,\,\}$.
Hence if $|\text{Re} ( z_{k})| > m$, then $|{\widehat f}(z_{k})|
 = |{\widehat f_{\text{Im} ( z_{k})}}(\text{Re} ( z_{k}))|
< \epsilon$.   So $z_{k} \rightarrow \infty$ implies that
$\eta_{z_{k}}(f) = {\widehat f}({z_{k}})  \rightarrow 0
= \eta_{\infty}(f)$.)
It follows as in the proof of Lemma 7.14 of \cite{\Ph} that the
compactly generated topology on $\text{Max}(A^{+})$ is the
topology it gets from the identification with $\varinjlim K_{n} $
(not the one point compactification of $\C$).
\par
By \cite{\Ph, Theorem 7.15}, $RK_{*}(A) = RK_{*}(C_{A})$, where $C_{A}$
is the closed subspace of $C(\text{Max}(A^{+}))$ of functions which
vanish at $\infty$.
It is easily seen, using the previous paragraph, that the restriction
map $C_{A} \rightarrow C_0 (\R)$ is a homotopy equivalence of Fr\'echet
algebras.  It follows that $A \rightarrow C_0 (\R)$ is an isomorphism
on $K$-theory.  Taking crossed products
by $\R$ and using the Thom isomorphism, we find that
$RK_{*} (\Cal S (\R, A, \alpha)) \rightarrow K_{*} (\Cal K)$ is an
isomorphism.
\par
We remark that the Gelfand transform maps any element of $A$ to a
Schwartz function on $\R$ which extends to a holomorphic
function on $\C$.  It follows that the
invertible elements of $A^{+}$ are precisely the nonzero
scalar multiples of the identity. Therefore any non-constant element in
$A^{+}$ has spectrum $\C$, and $A$ cannot be spectral invariant in any
Banach algebra.
\vskip\baselineskip

\enddocument